\pdfoutput=1

\documentclass[oneside]{article}
\usepackage{hyperref} % (find-es "tex" "hyperref")
  \usepackage[utf8]{inputenc}
\DeclareUnicodeCharacter{00B0}{^\circ}             % °
\DeclareUnicodeCharacter{00B9}{^{-1}}              % ¹
\DeclareUnicodeCharacter{00B2}{^2}                 % ²
\DeclareUnicodeCharacter{00B3}{^3}                 % ³
\DeclareUnicodeCharacter{00B1}{\pm}                % ±
\DeclareUnicodeCharacter{00F7}{\div}               % ÷
\DeclareUnicodeCharacter{00B7}{\cdot}              % ·
\DeclareUnicodeCharacter{00D7}{\times}             % ×
\DeclareUnicodeCharacter{00AC}{\neg}               % ¬
\DeclareUnicodeCharacter{00A7}{\S}                 % §
\DeclareUnicodeCharacter{00A1}{\text{\textexclamdown}} % ¡
\DeclareUnicodeCharacter{03BC}{\mu}                % μ
\DeclareUnicodeCharacter{03BD}{\nu}                % ν
\DeclareUnicodeCharacter{03C1}{\rho}               % ρ
\DeclareUnicodeCharacter{03C3}{\sigma}             % σ
\DeclareUnicodeCharacter{03C4}{\tau}               % τ
\DeclareUnicodeCharacter{03C6}{\phi}               % φ
\DeclareUnicodeCharacter{03C7}{\chi}               % χ
\DeclareUnicodeCharacter{03C8}{\psi}               % ψ
\DeclareUnicodeCharacter{0394}{\Delta}             % Δ
\DeclareUnicodeCharacter{0393}{\Gamma}             % Γ
\DeclareUnicodeCharacter{0398}{\Theta}             % Θ
\DeclareUnicodeCharacter{03A0}{\Pi}                % Π
\DeclareUnicodeCharacter{03A3}{\Sigma}             % Σ
\DeclareUnicodeCharacter{03A6}{\Phi}               % Φ
\DeclareUnicodeCharacter{03A8}{\Psi}               % Ψ
\DeclareUnicodeCharacter{03A9}{\Omega}             % Ω
\DeclareUnicodeCharacter{03B1}{\alpha}             % α
\DeclareUnicodeCharacter{03B2}{\beta}              % β
\DeclareUnicodeCharacter{03B3}{\gamma}             % γ
\DeclareUnicodeCharacter{03B4}{\delta}             % δ
\DeclareUnicodeCharacter{03B5}{\epsilon}           % ε
\DeclareUnicodeCharacter{03B7}{\eta}               % η
\DeclareUnicodeCharacter{03B8}{\theta}             % θ
\DeclareUnicodeCharacter{03B9}{\iota}              % ι
\DeclareUnicodeCharacter{03BB}{\lambda}            % λ
\DeclareUnicodeCharacter{03C0}{\pi}                % π
\DeclareUnicodeCharacter{03C9}{\omega}             % ω
\DeclareUnicodeCharacter{2022}{\bullet}            % •
\DeclareUnicodeCharacter{214B}{\bindnasrepma}      % ⅋
\DeclareUnicodeCharacter{2190}{\ot}                % ←
\DeclareUnicodeCharacter{2191}{\upto}              % ↑
\DeclareUnicodeCharacter{2192}{\to}                % →
\DeclareUnicodeCharacter{2193}{\dnto}              % ↓
\DeclareUnicodeCharacter{2194}{\bij}               % ↔
\DeclareUnicodeCharacter{2195}{\updownarrow}       % ↕
\DeclareUnicodeCharacter{2196}{\nwarrow}           % ↖
\DeclareUnicodeCharacter{2197}{\nearrow}           % ↗
\DeclareUnicodeCharacter{2198}{\searrow}           % ↘
\DeclareUnicodeCharacter{2199}{\swarrow}           % ↙
\DeclareUnicodeCharacter{21A3}{\epito}             % ↣
\DeclareUnicodeCharacter{21A3}{\twoheadrightarrow} % ↣
\DeclareUnicodeCharacter{21A4}{\mapsfrom}          % ↤
\DeclareUnicodeCharacter{21A6}{\mapsto}            % ↦
\DeclareUnicodeCharacter{21D0}{\Leftarrow}         % ⇐
\DeclareUnicodeCharacter{21D2}{\funto}             % ⇒
\DeclareUnicodeCharacter{2200}{\forall}            % ∀
\DeclareUnicodeCharacter{2203}{\exists}            % ∃
\DeclareUnicodeCharacter{2205}{\emptyset}          % ∅
\DeclareUnicodeCharacter{2208}{\in}                % ∈
\DeclareUnicodeCharacter{2216}{\backslash}         % ∖
\DeclareUnicodeCharacter{2218}{\circ}              % ∘
\DeclareUnicodeCharacter{221A}{\sqrt}              % √
\DeclareUnicodeCharacter{221E}{\infty}             % ∞
\DeclareUnicodeCharacter{2227}{\land}              % ∧
\DeclareUnicodeCharacter{2228}{\lor}               % ∨
\DeclareUnicodeCharacter{2229}{\cap}               % ∩
\DeclareUnicodeCharacter{2229}{\cap}               % ∩
\DeclareUnicodeCharacter{222A}{\cup}               % ∪
\DeclareUnicodeCharacter{2243}{\simeq}             % ≃
\DeclareUnicodeCharacter{2245}{\cong}              % ≅
\DeclareUnicodeCharacter{2260}{\neq}               % ≠
\DeclareUnicodeCharacter{2264}{\le}                % ≤
\DeclareUnicodeCharacter{2265}{\ge}                % ≥
\DeclareUnicodeCharacter{2282}{\subset}            % ⊂
\DeclareUnicodeCharacter{2283}{\supset}            % ⊃
\DeclareUnicodeCharacter{2286}{\subseteq}          % ⊆
\DeclareUnicodeCharacter{2287}{\supseteq}          % ⊇
\DeclareUnicodeCharacter{2293}{\sqcap}             % ⊓
\DeclareUnicodeCharacter{2294}{\sqcup}             % ⊔
\DeclareUnicodeCharacter{22A2}{\vdash}             % ⊢
\DeclareUnicodeCharacter{22A3}{\dashv}             % ⊣
\DeclareUnicodeCharacter{22A4}{\top}               % ⊤
\DeclareUnicodeCharacter{22A5}{\bot}               % ⊥
\DeclareUnicodeCharacter{22C0}{\bigwedge}          % ⋀
\DeclareUnicodeCharacter{22C1}{\bigvee}            % ⋁
\DeclareUnicodeCharacter{22C2}{\bigcap}            % ⋂
\DeclareUnicodeCharacter{22C3}{\bigcup}            % ⋃
\DeclareUnicodeCharacter{266D}{\flat}              % ♭
\DeclareUnicodeCharacter{266E}{\natural}           % ♮
\DeclareUnicodeCharacter{266F}{\sharp}             % ♯
\DeclareUnicodeCharacter{2806}{{:}}                % ⠆
\DeclareUnicodeCharacter{27E6}{\llbracket}         % ⟦
\DeclareUnicodeCharacter{27E7}{\rrbracket}         % ⟧
\DeclareUnicodeCharacter{223C}{\sim}               % ∼
\DeclareUnicodeCharacter{2248}{\approx}            % ≈
\DeclareUnicodeCharacter{2261}{\equiv}             % ≡
\DeclareUnicodeCharacter{22C4}{\lozenge}           % ⋄
\DeclareUnicodeCharacter{22C5}{\Box}               % ⋅
\DeclareUnicodeCharacter{25FB}{\Box}               % ◻
\DeclareUnicodeCharacter{25A1}{\Box}               % □
\DeclareUnicodeCharacter{2026}{\ldots}             % …
\DeclareUnicodeCharacter{222B}{\int}               % ∫
\DeclareUnicodeCharacter{21C0}{\rightharpoonup}    % ⇀
\DeclareUnicodeCharacter{2329}{\langle}            % 〈
\DeclareUnicodeCharacter{232A}{\rangle}            % 〉
\DeclareUnicodeCharacter{2299}{\odot}              % ⊙
\DeclareUnicodeCharacter{2296}{\ominus}            % ⊖
\DeclareUnicodeCharacter{2295}{\oplus}             % ⊕
\DeclareUnicodeCharacter{2298}{\oslash}            % ⊘
\DeclareUnicodeCharacter{2297}{\otimes}            % ⊗
\DeclareUnicodeCharacter{2581}{\_}                 % ▁
\DeclareUnicodeCharacter{1D41B}{\mathbf}           % 𝐛
\DeclareUnicodeCharacter{1D422}{\textsl}           % 𝐢
\DeclareUnicodeCharacter{1D42B}{\mathrm}           % 𝐫
\DeclareUnicodeCharacter{1D42D}{\text}             % 𝐭
\DeclareUnicodeCharacter{1D42C}{\mathsf}           % 𝐬

  \DeclareUnicodeCharacter{00B9}{^{*}}              % ¹
  \DeclareUnicodeCharacter{00B2}{^{**}}             % ²
  \DeclareUnicodeCharacter{00B3}{^{***}}            % ³
\usepackage{amsmath}
\usepackage{amsfonts}
\usepackage{amssymb}
\usepackage{stmaryrd}
\usepackage{pict2e}
\usepackage[x11names,svgnames]{xcolor} % (find-es "tex" "xcolor")
\usepackage{proof}   % For derivation trees ("%:" lines)
\input diagxy        % For 2D diagrams ("%D" lines)
\usepackage{edrx15}               % (find-LATEX "edrx15.sty")
\def\bsm#1{\left[\begin{smallmatrix}#1\end{smallmatrix}\right]}

\def\bsm#1{\left [\begin{smallmatrix}#1\end{smallmatrix}\right ]}

\def\subst#1{\left[\sm{#1}\right]}

\def\beginpicture (#1,#2)(#3,#4){\expr{beginpicture(v(#1,#2),v(#3,#4))}}
\def\beginpictureb(#1,#2)(#3,#4)#5{\expr{beginpicture(v(#1,#2),v(#3,#4),#5)}}

\def\pictaxes{\expr{pictaxes()}}
\def\pictaxes{{\linethickness{0.6pt}\expr{pictaxes()}}}
\def\pictaxes{{\linethickness{0.5pt}\expr{pictaxes()}}}
\def\pictgrid{{\color{GrayPale}\expr{pictgrid()}}}
\def\pictgrid{{\color{GrayPale}\linethickness{0.3pt}\expr{pictgrid()}}}
\def\pictdots#1{\expr{pictdots("#1")}}
\def\pictpiecewise#1{\expr{pictpiecewise("#1")}}
\def\pictpiecewise#1{{\linethickness{1pt}\expr{pictpiecewise("#1")}}}
\def\picturedots(#1,#2)(#3,#4)#5{%
  \vcenter{\hbox{%
  \beginpicture(#1,#2)(#3,#4)%
  \pictaxes%
  \pictdots{#5}%
  \end{picture}%
  }}%
}
\def\picturepiecewise(#1,#2)(#3,#4)#5{%
  \vcenter{\hbox{%
  \beginpicture(#1,#2)(#3,#4)%
  \pictgrid%
  \pictaxes%
  \pictpiecewise{#5}%
  \end{picture}%
  }}%
}

\def\Vector(#1,#2)(#3,#4){\expr{pict2evector(#1, #2, #3, #4)}}

\def\verteq{
 \vcenter{\hbox{%
 \unitlength=1pt%
 \linethickness{0.4pt}%
 \beginpicture(0,0)(5,6)
   \roundcap
   \Line(1,0)(1,6)
   \Line(3,0)(3,6)
 \end{picture}%
 }}%
}

\usepackage[backend=biber,
   style=alphabetic]{biblatex} % (find-es "tex" "biber")
\usepackage[a4paper]{geometry}             % (find-es "tex" "geometry")
\begin{document}

% This file: (find-LATEX "2017planar-has-defs.tex")
% «.picturedots»	(to "picturedots")
% «.defs»		(to "defs")
% «.squigbij»		(to "squigbij")
% «.squigbijtest»	(to "squigbijtest")
% «.defzha-and-deftcg»	(to "defzha-and-deftcg")
% «.defpido»		(to "defpido")
% «.defpicturedots»	(to "defpicturedots")
% «.picturedotsdef»	(to "picturedotsdef")
% «.defub»		(to "defub")

%        _      _                      _       _       
%  _ __ (_) ___| |_ _   _ _ __ ___  __| | ___ | |_ ___ 
% | '_ \| |/ __| __| | | | '__/ _ \/ _` |/ _ \| __/ __|
% | |_) | | (__| |_| |_| | | |  __/ (_| | (_) | |_\__ \
% | .__/|_|\___|\__|\__,_|_|  \___|\__,_|\___/ \__|___/
% |_|                                                  
%
% «picturedots» (to ".picturedots")
% (find-LATEX "edrxpict.lua" "beginpicture")
% (find-LATEX "edrxpict.lua" "pictaxes")
% (find-LATEX "edrxpict.lua" "pictdots")
% (find-LATEX "2016-2-GA-algebra.tex" "picturedots")
% (find-LATEX "2016-2-GA-algebra.tex" "comprehension-gab")
% (gaap 5)
%
\def\beginpicture(#1,#2)(#3,#4){\expr{beginpicture(v(#1,#2),v(#3,#4))}}
\def\pictaxes{\expr{pictaxes()}}
\def\pictdots#1{\expr{pictdots("#1")}}
\def\picturedotsa(#1,#2)(#3,#4)#5{%
  \vcenter{\hbox{%
  \beginpicture(#1,#2)(#3,#4)%
  \pictaxes%
  \pictdots{#5}%
  \end{picture}%
  }}%
}
\def\picturedots(#1,#2)(#3,#4)#5{%
  \vcenter{\hbox{%
  \beginpicture(#1,#2)(#3,#4)%
  %\pictaxes%
  \pictdots{#5}%
  \end{picture}%
  }}%
}

%      _       __     
%   __| | ___ / _|___ 
%  / _` |/ _ \ |_/ __|
% | (_| |  __/  _\__ \
%  \__,_|\___|_| |___/
%                     
% «defs» (to ".defs")

\def\sa{\rightsquigarrow}
\def\BPM{\mathsf{BPM}}
\def\WPM{\mathsf{WPM}}
\def\ZHAG{\mathsf{ZHAG}}

\def\LR{\mathbb{LR}}

\def\catTwo{\mathbf{2}}
\def\calS{\mathcal{S}}
\def\calI{\mathcal{I}}
\def\calK{\mathcal{K}}
\def\calV{\mathcal{V}}

\def\und#1#2{\underbrace{#1}_{#2}}

\def\subst#1{\left[\begin{array}{rcl}#1\end{array}\right]}
\def\subst{\bsm}

% (find-LATEXfile "2015planar-has.tex" "\\def\\Mop")

\def\MP  {\mathsf{MP}}
\def\J   {\mathsf{J}}
\def\Mo  {\mathsf{Mo}}
\def\Mop {\mathsf{Mop}}
\def\Sand{\mathsf{Sand}}
\def\ECa {\mathsf{EC}{\&}}
\def\ECv {\mathsf{EC}{∨}}
\def\ECS {\mathsf{ECS}}
\def\pdiag#1{\left(\diag{#1}\right)}
\def\ltor#1#2{#1\_{\to}\_#2}
\def\lotr#1#2{#1\_{\ot}\_#2}
\def\Int{{\operatorname{int}}}
\def\Int{{\operatorname{\mathsf{int}}}}
\def\coInt{{\operatorname{\mathsf{coint}}}}
\def\LC {\mathsf{LC}}
\def\RC {\mathsf{RC}}
\def\TCG{\mathsf{2CG}}
\def\pile{\mathsf{pile}}
\def\ltor#1#2{#1\_{\to}\_#2}
\def\lotr#1#2{#1\_{\ot}\_#2}
\def\ltol#1#2{#1\_{\to}#2\_}
\def\rtor#1#2{\_#1{\to}\_#2}
\def\NoLcuts{\mathsf{No}λ\mathsf{cuts}}
\def\NoYcuts{\mathsf{NoYcuts}}
\def\astarcube{{\&}^*\mathsf{Cube}}
\def\ostarcube{{∨}^*\mathsf{Cube}}
\def\istarcube{{→}^*\mathsf{Cube}}
\def\acz{{\&}^*\mathsf{C}_0}
\def\ocz{{∨}^*\mathsf{C}_0}
\def\icz{{→}^*\mathsf{C}_0}
\def\astarcuben{{\&}^*\mathsf{Cube}_\mathsf{n}}
\def\ostarcuben{{∨}^*\mathsf{Cube}_\mathsf{n}}
\def\istarcuben{{→}^*\mathsf{Cube}_\mathsf{n}}
\def\astarcubev{{\&}^*\mathsf{Cube}_\mathsf{v}}
\def\ostarcubev{{∨}^*\mathsf{Cube}_\mathsf{v}}
\def\istarcubev{{→}^*\mathsf{Cube}_\mathsf{v}}
%
%\catcode`∧=13 \def∧{\mathop{\&}}

\def\biggest {\mathsf{biggest}}
\def\smallest{\mathsf{smallest}}
\def\Cuts    {\mathsf{Cuts}}

\def\myresizebox#1{%
  \noindent\hbox to \textwidth{\hss
    \resizebox{1.0\textwidth}{!}{#1}%
    \hss}
  }

\def\LR  {\mathbb{LR}}
\def\Taut{\mathsf{Taut}}
\def\IPL {\mathrm{IPL}}
\def\CPL {\mathrm{CPL}}
\def\ZHAL{\mathrm{ZHAL}}

% «squigbij» (to ".squigbij")
% (ph2 "question-marks")
% (find-es "tex" "pict2e-squigbij")

\def\squigbij{\newsquigbij}
\def\oldsquigbij{\;\; \diagxyto/<~>/<300> \;\;}

\def\newsquigbij{\;\; \squigbijbody \;\;}
\def\squigbijy{-1.2}
\def\squigbijbody{\squigbijbodywithparams{1.5pt}{0.3pt}{1.0}}
\def\squigbijtriangle(#1,#2)#3{\polygon*(#1,0)(#2,#3)(#2,-#3)}
\def\squigbijbodywithparams#1#2#3{{%
  \unitlength=#1
  \linethickness{#2}
  % \beginpicture(-5,-1)(17,1)%
  % \begin{picture}(22.4,2.4)(-5.2,-1.2)%
  \begin{picture}(22.4,2.4)(-5.2,\squigbijy)%
    \polyline(-3,0)(0,0)(1,1)(3,-1)(5,1)(7,-1)(9,1)(11,-1)(12,0)(14,0)
    \squigbijtriangle(-5,-2){#3}
    \squigbijtriangle(17,14){#3}
  \end{picture}%
  }}

% «squigbijtest»  (to ".squigbijtest")
\def\squigbijtest#1{
  \def\squigbijy{#1}
  \par #1: $A \squigbij B$
  }
% Use something like this to find a good value for \squigbijy.
%   \squigbijtest{-4.0}
%   \squigbijtest{-3.0}
%   \squigbijtest{-2.0}
%   \squigbijtest{-1.0}
%   \squigbijtest{0.0}
%   \squigbijtest{1.0}
%   \squigbijtest{2.0}
%   \squigbijtest{3.0}
%   \squigbijtest{4.0}
% These values work well for me:
%   \squigbijtest {-2.5}   % for 12pt
%   \def\squigbijy{-2.5}   % for 12pt
%   \squigbijtest {1.2}    % for 10pt
%   \def\squigbijy{1.2}    % for 10pt

% «defzha-and-deftcg»  (to ".defzha-and-deftcg")
% (find-es "dednat" "defzha-and-deftcg")
\def\defzha#1#2{\expandafter\def\csname zha-#1\endcsname{#2}}
\def\ifzhaundefined#1{\expandafter\ifx\csname zha-#1\endcsname\relax}
\def\zha#1{\ifzhaundefined{#1}
    \errmessage{UNDEFINED ZHA: #1}
  \else
    \csname zha-#1\endcsname
  \fi
}
\def\deftcg#1#2{\expandafter\def\csname tcg-#1\endcsname{#2}}
\def\iftcgundefined#1{\expandafter\ifx\csname tcg-#1\endcsname\relax}
\def\tcg#1{\iftcgundefined{#1}
    \errmessage{UNDEFINED TCG: #1}
  \else
    \csname tcg-#1\endcsname
  \fi
}

% «defpido»  (to ".defpido")
% (find-LATEX "edrxpict.lua" "defpictdots")
% (find-LATEX "2019oxford-abs.tex" "defpictdots")
\def\defpido#1#2{\expandafter\def\csname pido-#1\endcsname{#2}}
\def\ifpidoundefined#1{\expandafter\ifx\csname pido-#1\endcsname\relax}
\def\pido#1{\ifpidoundefined{#1}
    \errmessage{UNDEFINED PIDO: #1}
  \else
    \csname pido-#1\endcsname
  \fi
}

% Not used?
\def\defpicturedots #1(#2,#3)(#4,#5)#6{%
    \directlua{defpictdots(nil, "#1", #2,#3, #4,#5,nil, "#6")}
  }
\def\defpicturedotsa#1(#2,#3)(#4,#5)#6{%
    \directlua{defpictdots("axes", "#1", #2,#3, #4,#5,nil, "#6")}
  }

% «picturedotsdef»  (to ".picturedotsdef")
% (find-LATEX "edrxpict.lua" "defpictdots")
% This lets me rewrite the first two lines below as the other two lines...
% \picturedotsa        (-2,0)(2,5){ 0,0 1,1 2,2  -1,1 0,2 1,3  -2,2 -1,3 0,4   -1,5 }
% \picturedots         (-2,0)(2,5){ 0,0 1,1 2,2  -1,1 0,2 1,3  -2,2 -1,3 0,4   -1,5 }
% \picturedotsadef{foo}(-2,0)(2,5){ 0,0 1,1 2,2  -1,1 0,2 1,3  -2,2 -1,3 0,4   -1,5 }
% \picturedotsdef {bar}(-2,0)(2,5){ 0,0 1,1 2,2  -1,1 0,2 1,3  -2,2 -1,3 0,4   -1,5 }
%
\def\picturedotsadef#1(#2,#3)(#4,#5)#6{
  \directlua{ defpictdots("axes", "#1", #2,#3, #4,#5, nil, "#6") }
  \pido{#1}
  }
\def\picturedotsdef #1(#2,#3)(#4,#5)#6{
  \directlua{ defpictdots(nil,    "#1", #2,#3, #4,#5, nil, "#6") }
  \pido{#1}
  }

% «defub»  (to ".defub")
% Example: (find-LATEX "2017planar-has-1.tex" "prop-calc-ZHA")
%          (find-LATEX "2017planar-has-1.tex" "prop-calc-ZHA" "defub")
\def\defub#1#2{\expandafter\def\csname ub-#1\endcsname{#2}}
\def\ifubundefined#1{\expandafter\ifx\csname ub-#1\endcsname\relax}
\def\ub#1{\ifubundefined{#1}
    \errmessage{UNDEFINED UB: #1}
  \else
    \csname ub-#1\endcsname
  \fi
}

% Local Variables:
% coding: utf-8-unix
% ee-anchor-format: "«%s»"
% End:
 % (find-LATEX "2017planar-has-defs.tex")
% (find-LATEX "2019J-ops-defs.tex")
% (defun c () (interactive) (find-LATEXsh "lualatex -record 2019J-ops-defs.tex" :end))
% (defun d () (interactive) (find-pdf-page "~/LATEX/2019J-ops-defs.pdf"))
% (defun e () (interactive) (find-LATEX "2019J-ops-defs.tex"))
% (defun u () (interactive) (find-latex-upload-links "2019J-ops-defs"))
% (find-pdf-page   "~/LATEX/2019J-ops-defs.pdf")
% (find-sh0 "cp -v  ~/LATEX/2019J-ops-defs.pdf /tmp/")
% (find-sh0 "cp -v  ~/LATEX/2019J-ops-defs.pdf /tmp/pen/")
%   file:///home/edrx/LATEX/2019J-ops-defs.pdf
%               file:///tmp/2019J-ops-defs.pdf
%           file:///tmp/pen/2019J-ops-defs.pdf
% http://angg.twu.net/LATEX/2019J-ops-defs.pdf
% (find-LATEX "2019.mk")

% «.defs»		(to "defs")
% «.chars»		(to "chars")

%  ____        __     
% |  _ \  ___ / _|___ 
% | | | |/ _ \ |_/ __|
% | |_| |  __/  _\__ \
% |____/ \___|_| |___/
%                     
% «defs»  (to ".defs")

\def\eqP{\underset{P}{\sim}}
\def\eqJ{\underset{J}{\sim}}
\def\eqP{\underset{\scriptscriptstyle P}{\sim}}
\def\eqJ{\underset{\scriptscriptstyle J}{\sim}}
\def\eqP{\sim_P}
\def\eqJ{\sim_J}
\def\eqL{\sim_L}
\def\eqR{\sim_R}
\def\eqS{\sim_S}
\def\eqF{\sim_F}
\def\eqQ{\sim_Q}
\def\eqQp{\sim_{Q'}}

\def\ECa {\mathsf{EC}{∧}}

\def\ECube{\mathsf{ECube}}   \def\ecube{\mathsf{ecube}}
\def\OCube{\mathsf{OCube}}   \def\ocube{\mathsf{ocube}}
\def\FCube{\mathsf{FCube}}   \def\fcube{\mathsf{fcube}}
\def\SCube{\mathsf{SCube}}   \def\fcube{\mathsf{fcube}}
\def\VCube{\mathsf{VCube}}   \def\vcube{\mathsf{vcube}}
\def\Exprs{\mathsf{Exprs}}
\def\Thms{\mathsf{Thms}}
\def\thms{\mathsf{thms}}
\def\vthms{\mathsf{vthms}}
\def\NClasses{\mathsf{NClasses}}
\def\nclasses{\mathsf{nclasses}}
\def\ZHAstar{ZHA${}^*$}
\def\sfE{\mathsf{E}}
\def\sfV{\mathsf{V}}

\def\oand{\varowedge}
\def\oor {\varovee}
\def\oimp{\mathbin{\buildoimp{\ominus}{\to}}}
\def\buildoimp#1#2{\rlap{$#1$}\hbox{$#2$}}

\def\zs{{0\ldots7}}

% «chars»  (to ".chars")
% (find-LATEX "2019oxford-chars.tex")
% \ifluatex
%   \catcode`¹=13 \def¹{^{*}}
%   \catcode`²=13 \def²{^{**}}
%   \catcode`³=13 \def³{^{***}}
% \else
%   \DeclareUnicodeCharacter{00B9}{^{*}}              % ¹
%   \DeclareUnicodeCharacter{00B2}{^{**}}             % ²
%   \DeclareUnicodeCharacter{00B3}{^{***}}            % ³
% \fi

% (find-eev "eev-math-glyphs.el" "faces")
% (eev-glyphs-set-face 'eev-glyph-face-math  "RoyalBlue2" "gray20")
% (eev-glyphs-set-face 'eev-glyph-face-math2 "DarkOrange" "gray20")
% (eev-glyphs-set-face 'eev-glyph-face-math2 "DarkOrange" "gray25")
% (eev-glyphs-set-face 'eev-glyph-face-math2 "DarkOrchid2" "gray25")
% (eev-glyphs-set-face 'eev-glyph-face-math2 "DarkOrchid2" "DarkOliveGreen")
% (eev-glyphs-set-face 'eev-glyph-face-math2 "DarkOrchid2" "blue4")
% (eepitch-set-glyph0 ?² ?² 'eev-glyph-face-math)
% (eepitch-set-glyph0 ?² ?² 'eev-glyph-face-math2)
% (progn (find-epalette) (ee-goto-position "RoyalBlue2"))
% (find-LATEX "2017planar-has-defs.tex" "defpido")
% ¹²

% Local Variables:
% coding: utf-8-unix
% ee-tla: "jod"
% End:
      % (find-LATEX "2019J-ops-defs.tex")

%\ifluatex
%  \catcode`\^^J=10
%  \directlua{dofile "dednat6load.lua"}
%\else
  \input\jobname.dnt
  \def\pu{}
  \def\directlua#1{}
%\fi

%L dofile "edrxtikz.lua"  -- (find-LATEX "edrxtikz.lua")
%L dofile "edrxpict.lua"  -- (find-LATEX "edrxpict.lua")
\pu

\def\ovl{\overline}

%  _____ _ _   _      
% |_   _(_) |_| | ___ 
%   | | | | __| |/ _ \
%   | | | | |_| |  __/
%   |_| |_|\__|_|\___|
%                     
% (find-LATEX "idarct/idarct-preprint.tex")
% «title» (to ".title")
% (jopp 1 "title")
% (jop    "title")

\title{Planar Heyting Algebras for Children 2: Local Operators, J-Operators, and Slashings}

\author{Eduardo Ochs}

\maketitle

%     _    _         _                  _   
%    / \  | |__  ___| |_ _ __ __ _  ___| |_ 
%   / _ \ | '_ \/ __| __| '__/ _` |/ __| __|
%  / ___ \| |_) \__ \ |_| | | (_| | (__| |_ 
% /_/   \_\_.__/|___/\__|_|  \__,_|\___|\__|
%                                           
% «abstract» (to ".abstract")
% (jopp 1 "abstract")
% (jop    "abstract")

\begin{abstract}

Choose a topos $\calE$. There are several different ``notions of
sheafness'' on $\calE$. How do we visualize them?

Let's refer to the classifier object of $\calE$ as $Ω$, and to its
Heyting Algebra of truth-values, $\Sub(1_\calE)$, as $H$; we will
sometimes call $H$ the ``logic'' of the topos. There is a well-known
way of representing notions of sheafness as morphisms $j:Ω→Ω$, but
these `$j$'s yield big diagrams when we draw them explicitly; here we
will see a way to represent these `$j$'s as maps $J:H→H$ in a way that
is much more manageable.

In the previous paper of this series --- called \cite{OchsPH1} from
here on --- we showed how certain toy models of Heyting Algebras,
called ``ZHAs'', can be used to develop visual intuition for how
Heyting Algebras and Intuitionistic Propositional Logic work; here we
will extend that to sheaves. The full idea is this: {\sl notions of
  sheafness} correspond to {\sl local operators} and vice-versa; {\sl
  local operators} correspond to {\sl J-operators} and vice-versa; if
our Heyting Algebra $H$ is a ZHA then {\sl J-operators} correspond to
{\sl slashings} on $H$, and vice-versa; {\sl slashings} on $H$
correspond to {\sl ``sets of question marks''} and vice-versa, and
each set of question marks induces a notion of {\sl erasing and
  reconstructing}, which induces a sheaf. Also, every ZHA $H$
corresponds to an {(acyclic) 2-column graph}, and vice-versa, and for
any two-column graph $(P,A)$ the logic of the topos $\Set^{(P,A)}$ is
exactly the ZHA $H$ associated to $(P,A)$.

The introduction of \cite{OchsPH1} discusses two different senses in
which a mathematical text can be ``for children''. The first sense
involves some precise metamathetical tools for transfering knowledge
back and forth between a general case ``for adults'' and a toy model
``for children''; the second sense is simply that the text's
presentation has few prerequisites and never becomes too abstract.
Here we will use the second sense: everything here, except for the
last section, should be accessible to students who have taken a course
on Discrete Mathematics and read \cite{OchsPH1}. This means that
categories, toposes, sheaves and the maps $j:Ω→Ω$ only appear in the
last section, and before that we deal only with the J-operators
$J:H→H$, how they correspond to slashings and sets of question marks,
and how they form an algebra.

\end{abstract}

% (find-books "__cats/__cats.el" "bell")
% (find-belltpage (+ 14 163) "modality")
% (find-belltpage (+ 14 164) "universal closure operation")
%
% (find-books "__cats/__cats.el" "johnstone-elephant")
% (find-elephantpage (+ 17 195) "A4.4 Local Operators")
%
% (find-books "__cats/__cats.el" "fourman-scott")
% (find-slnm0753page (+ 16 324) "J-operators")

% TODO: add this (and an abstract):
% (ph2p 4 "piccs-and-slashings")
% (ph2    "piccs-and-slashings")

% \end{document}

%  ____            _       
% |  _ \ __ _ _ __| |_ ___ 
% | |_) / _` | '__| __/ __|
% |  __/ (_| | |  | |_\__ \
% |_|   \__,_|_|   \__|___/
%                          

% «parts»  (to ".parts")

\newpage
% (find-LATEX "2019J-ops-slashings.tex")
% (defun c () (interactive) (find-LATEXsh "lualatex -record 2019J-ops-slashings.tex" :end))
% (defun d () (interactive) (find-pdf-page "~/LATEX/2019J-ops-slashings.pdf"))
% (defun e () (interactive) (find-LATEX "2019J-ops-slashings.tex"))
% (defun u () (interactive) (find-latex-upload-links "2019J-ops-slashings"))
% (find-pdf-page   "~/LATEX/2019J-ops-slashings.pdf")
% (find-sh0 "cp -v  ~/LATEX/2019J-ops-slashings.pdf /tmp/")
% (find-sh0 "cp -v  ~/LATEX/2019J-ops-slashings.pdf /tmp/pen/")
%   file:///home/edrx/LATEX/2019J-ops-slashings.pdf
%               file:///tmp/2019J-ops-slashings.pdf
%           file:///tmp/pen/2019J-ops-slashings.pdf
% http://angg.twu.net/LATEX/2019J-ops-slashings.pdf
% (find-LATEX "2019.mk")

% «.basic-definitions»		(to "basic-definitions")
% «.qms-and-slashings»		(to "qms-and-slashings")
% «.piccs-and-slashings»	(to "piccs-and-slashings")
% «.slash-ops»			(to "slash-ops")
% «.converting»			(to "converting")

\directlua{tf_push("2019J-ops-slashings.tex")}

%  ____            _            _       __     
% | __ )  __ _ ___(_) ___    __| | ___ / _|___ 
% |  _ \ / _` / __| |/ __|  / _` |/ _ \ |_/ __|
% | |_) | (_| \__ \ | (__  | (_| |  __/  _\__ \
% |____/ \__,_|___/_|\___|  \__,_|\___|_| |___/
%                                              
% «basic-definitions»  (to ".basic-definitions")
% (jopp 2 "basic-definitions")
% (jos    "basic-definitions")

%\section{Question marks and slashings}
\section{Basic definitions}
\label  {basic-definitions}

One of the main constructions of \cite{OchsPH1} is a correspondence
between 2-column graphs (``2CGs'') and Planar Heyting Algebras
(``ZHAs''), as in this example:
%
% (oxap 7 "fig:2CGs-ZHA")
% (oxa    "fig:2CGs-ZHA")
% (oxa    "fig:2CGs-ZHA" "\\tcg{(P,A)}")
% 
%L tdims = TCGDims {qrh=5, q=15, crh=12, h=60, v=25, crv=7}   -- with v arrows
%L tspec_PA  = TCGSpec.new("46; 11 22 34 45, 25")
%L tspec_PAQ = TCGSpec.new("46; 11 22 34 45, 25", ".???", "???.?.")
%L tspec_PA :mp  ({zdef="O_A(P)"})  :addlrs():print()            :output()
%L tspec_PAQ:mp  ({zdef="O_A(P),J"}):addlrs():print()            :output()
%L tspec_PA :tcgq({tdef="(P,A)",   meta="1pt p"}, "lr q h v ap") :output()
%L tspec_PAQ:tcgq({tdef="(P,A),Q", meta="1pt p"}, "lr q h v ap") :output()
%L
%L tspec_PAC = TCGSpec.new("46; 11 22 34 45, 25", "?...", "???...")
%L tspec_PAC:mp  ({zdef="closed-op"}) :addlrs():print()            :output()
%L tspec_PAC:tcgq({tdef="closed-op", meta="1pt p"}, "lr q h v ap") :output()
%L 
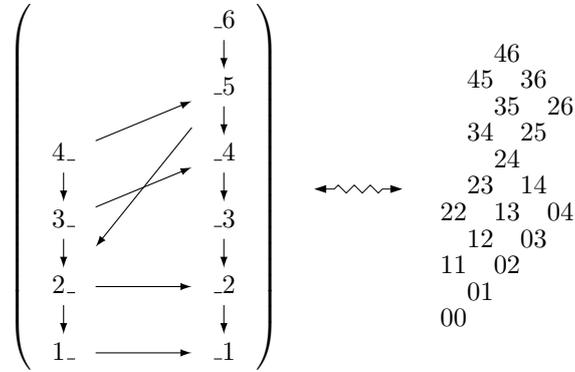
\pu
$$\tcg{(P,A)} \;\; \squigbij \;\;\; \zha{O_A(P)}$$

The arrows in the 2CG $(P,A)$ (mnemonic: ``points'' and ``arrows'')
are interpreted as conditions that subsets of $P$ must obey to the
open: for example, the arrow $(4\_,\_5)∈A$ means that if an open set
$U⊆A$ contains the point $4\_$ then it also has to contains $\_5$.
This generates an {\sl order topology} on $P$, that we denote by
$\Opens_A(P)$, and the ZHA $H$ at the right of the squiggly arrow in
the figure is this $\Opens_A(P)$ drawn in a very compacty way --- by
using the operation ``$\pile$'', and abbreviating it.

We write $\pile(ab)$ for the subset of $P$ formed by pile of $a$
elements at the left and a pile of $b$ elements at the right, as in:
$$25 ≡ \pile(25) = \{2▁,1▁, \;\; ▁1,▁2,▁3,▁4,▁5\},$$
The `$≡$' in ``$25 ≡ \pile(2,5)$'' means a change of notation --- it
means that sometimes `$ab$' will be an abbreviation for
``$\pile(ab)$''. With this abreviation it is easy to check that the
$H$ above is exactly the topology $\Opens_A(P)$. Note that, for
example, $21\not∈H$; this is because $\pile(21) = \{2▁,1▁, \;\; ▁1\}$,
and this set does not obey all the conditions associated to the arrows
in $A$: we have $(2\_,\_2)∈A$ but $2\_∈\pile(21)$ and
$\_2\not∈\pile(21)$.

Let's now introduce some new ideas.

%   ___                                  _       _     
%  / _ \ _ __ ___  ___    __ _ _ __   __| |  ___| |___ 
% | | | | '_ ` _ \/ __|  / _` | '_ \ / _` | / __| / __|
% | |_| | | | | | \__ \ | (_| | | | | (_| | \__ \ \__ \
%  \__\_\_| |_| |_|___/  \__,_|_| |_|\__,_| |___/_|___/
%                                                      
% «qms-and-slashings»  (to ".qms-and-slashings")
% (jopp 2 "qms-and-slashings")
% (jos    "qms-and-slashings")
\subsection{Question marks and slashings}
\label    {qms-and-slashings}

A {\sl set of question marks} on a 2CG $(P,A)$ is a subset $Q⊆P$. We
write a 2CG with question marks as $((P,A),Q)$, and we represent this
$Q$ graphically by writing a `?' close to each element of $P$ that
belongs to $Q$, as in the figure below. The intended meaning of these
question marks is that we want to forget the information on them and
then see which elements of $\Opens_A(P)$ become indistinguishable
after this forgetting: two elements $ab,cd∈H$ are {\sl
$Q$-equivalent}, written as $ab \eqQ cd$, iff $\pile(ab)∖Q =
\pile(cd)∖Q$. In the $((P,A),Q)$ of the figure below we have $23 \eqQ
13 \not\eqQ 14$.

A {\sl slashing} $S$ on a ZHA $H$ is a set of diagonal cuts on $H$
``that do not stop midway''. These cuts are interpreted as fences that
divide $H$ in separate regions, and two elements $ab,cd∈H$ are {\sl
$S$-equivalent}, written as $ab \eqS cd$, if they belong to the same
region. In the slashing at the right in the figure below we have $11
\eqS 23 \not\eqS 14$.
$$\tcg{(P,A),Q} \;\; \squigbij \;\;\; \zha{O_A(P),J}$$

In \cite{OchsPH1} we used the notation $(P,A) \; \squigbijbody \; H$
to say that $H$ is the ZHA associated to the 2CG $(P,A)$; this ``is
associated to'' was interpreted formally as $\Opens_A(P) = H$. We are
now extending this to $((P,A),Q) \; \squigbijbody \; (H,S)$ --- a 2CG
with question marks $((P,A),Q)$ is associated to the ZHA with slashing
$(H,S)$ when we have $\Opens_A(P) = H$ and the equivalence relations
$\eqQ,\eqS⊆H×H$ coincide. Note that the two `$\squigbijbody$'s are
both pronounced as ``is associated to'', but they have different
formal meanings.

%  ____  _                                _       _     _         
% |  _ \(_) ___ ___ ___    __ _ _ __   __| |  ___| |___| |__  ___ 
% | |_) | |/ __/ __/ __|  / _` | '_ \ / _` | / __| / __| '_ \/ __|
% |  __/| | (_| (__\__ \ | (_| | | | | (_| | \__ \ \__ \ | | \__ \
% |_|   |_|\___\___|___/  \__,_|_| |_|\__,_| |___/_|___/_| |_|___/
%                                                                 
% «piccs-and-slashings»  (to ".piccs-and-slashings")
% (jopp 3 "piccs-and-slashings")
% (jos    "piccs-and-slashings")
% (ph2p 4 "piccs-and-slashings")
% (ph2    "piccs-and-slashings")
\subsection{Piccs and slashings}
\label     {piccs-and-slashings}

A picc (``partition into contiguous classes'') of a ``discrete
interval'' $I=\{0,\ldots,n\}$ is a partition $P$ of $I$ that obeys
this condition (``picc-ness''):
$$∀a,b,c∈\{0,\ldots,n\}.\; (a<b<c ∧ a \eqP c) → (a \eqP b ∧ b \eqP c).$$
So $P = \{\{0\},\{1,2,3\},\{4,5\}\}$ is a picc of $\{0,\ldots,5\}$,
and $P' = \{\{0\},\{1,2,4,5\},\{3\}\}$ is a partition of
$\{0,\ldots,5\}$ that is not a picc.

A short notation for piccs is this:
$$0|123|45 \equiv \{\{0\},\{1,2,3\},\{4,5\}\}$$
we list all digits in the (discrete) interval in order, and we put
bars to indicate where we change from one equivalence class to
another.

\msk

We will represent a slashing $S$ formally as pairs of piccs, one for
the left digit and one for the right digit. Our notation for slashings
as pairs will be based on this figure:
%
%L -- (find-LATEX "dednat6/zhas.lua" "VCuts-tests")
%L local vc = VCuts.new({scale="7pt", def="VCuts"}, 4, 6)
%L vc:cutl(0)
%L vc:cutr(3):cutr(5)
%L vc:output()
%L
%L mp = mpnew({def="ZQuot"},      "12345RR4321"):addlrs():addcuts("c 4321/0 0123|45|6"):output()
%L mp = mpnew({def="ZQuotients"}, "1R2R3212RL1"):addlrs():addcuts("c 4321/0 0123|45|6"):output()
%L mp:print()
%
$$\pu
  \VCuts
  \qquad
  % \ZQuot
  % \qquad
  \ZQuotients
$$

The slashing $S$ that we are using in our examples will be represented
as:
$$\begin{array}{rcl}
  S &=& (L,R) \\
    &=& (\{\{0\},\{1,2,3,4\}\}, \, \{\{0,1,2,3\},\{4,5\},\{6\}\}) \\
    &=& (0|1234, 0123|45|6) \\
    &=& (4321/0,\, 0123∖45∖6) \\
  \end{array}
$$
We use `$/$'s and `$∖$'s instead of `$|$'s to remind us of the
direction of the cuts: the `$/$'s correspond to cuts that go northeast
and the `$∖$'s to cuts that go northwest.

We can now define the equivalence relation $\eqS$ formally: if
$S=(L,R)$ then $ab \eqS cd$ iff $a \eqL c$ and $c \eqR d$.

\msk

The expression ``$S=(L,R)$ is a slashing on $H$'' will mean: $H$ is a
ZHA, $L$ is a picc on $\{0,\ldots,l\}$, and $R$ a picc on
$\{0,\ldots,r\}$, where $lr$ is the top element of $H$. The domain of
the equivalence relation $\eqS$ will be considered to be $H$, not
$\{0,\ldots,l\} × \{0,\ldots,r\}$.

% Note that $\eqS$ is an equivalence relation on $H$, not on
% $\{0,\ldots,a\} × \{0,\ldots,b\}$.

%  ____  _           _                           
% / ___|| | __ _ ___| |__         ___  _ __  ___ 
% \___ \| |/ _` / __| '_ \ _____ / _ \| '_ \/ __|
%  ___) | | (_| \__ \ | | |_____| (_) | |_) \__ \
% |____/|_|\__,_|___/_| |_|      \___/| .__/|___/
%                                     |_|        
%
% «slash-ops»  (to ".slash-ops")
% (jopp 4 "slash-ops")
% (jos    "slash-ops")
\subsection{Slash-operators}
\label    {slash-ops}

When $S=(L,R)$ is a slashing on $H$ we will use the notations $[·]^L$,
$[·]^R$, $[·]^S$ for the equivalence classes of $L$, $R$, $S$ and the
notations $·^L$, $·^R$, $·^S$ for the highest element in those
equivalence class. In our example we have $[2]^L = \{1,2,3,4\}$,
$[2]^R = \{0,1,2,3\}$, $[22]^S = \{11,12,13,22,23\}$, $2^L = 4$, $2^R
= 4$, $2^S = 23$. Note that $[·]^S$ and $·^S$ depend on the ZHA.

A {\sl slash-operator} on a ZHA $H$ is a function $f:H→H$ that is
equal to some $·^S$.

\msk

Take any function $f:H→H$ on a ZHA. Let:
$$\begin{array}{rcl}
  S_0 &=& \setofst{(ab,f(ab))}{ab∈H} \\
  L_0 &=& \setofst{(a,c)}{(ab,cd)∈S_0} \\
  R_0 &=& \setofst{(b,d)}{(ab,cd)∈S_0} \\
  L   &=& {L_0}^* \\
  R   &=& {R_0}^* \\
  S   &=& (L,R) \\
  \end{array}
$$
The function $f$ is a slash-operator if and only if these $L$ and $R$
are piccs and $f = ·^S$.

%   ____                          _   _             
%  / ___|___  _ ____   _____ _ __| |_(_)_ __   __ _ 
% | |   / _ \| '_ \ \ / / _ \ '__| __| | '_ \ / _` |
% | |__| (_) | | | \ V /  __/ |  | |_| | | | | (_| |
%  \____\___/|_| |_|\_/ \___|_|   \__|_|_| |_|\__, |
%                                             |___/ 
%
% «converting»  (to ".converting")
% (jopp 4 "converting")
% (jos    "converting")

\subsection{From slashings to question marks and vice-versa}

Choose any path from the bottom element of the ZHA to its top element
that is made of one unit steps northwest or northeast --- for example,
this one:
$$(a_0b_0, a_1b_1, \ldots a_{10}b_{10}) = (00, 01, 02, 03, 04, 14, 24, 34, 35, 36, 46)$$

If we apply `$\pile$' to each element of that path we get a sequence
of sets,
$$(\pile(a_0b_0), \pile(a_1b_1), \ldots, \pile(a_{10}b_{10}))$$
that is actually a sequence of open sets in $\Opens_A(P)$ in which the
first set is $\pile(a_0b_0) = \pile(00) = ∅$, the last set is $P$, and
the difference between each set and the next one is exactly one
element --- for example:
$$\begin{array}{rcl}
  \pile(34)∖\pile(24) &=& \{3▁\} \\
  \pile(35)∖\pile(34) &=& \{▁5\} \\
  \end{array}
$$

{

\def\aibi   {a_ib_i}
\def\aiibii {a_{i+1}b_{i+1}}
\def\paibi  {\pile(\aibi)}
\def\paiibii{\pile(\aiibii)}

Note that we have two different cases: 1) the step from $\aibi$ to
$\aiibii$ is a movement {\sl northwest} in the ZHA, as in from 24 to
34; in this case $a_{i+1}b_{i+1} = (a_i+1)b_i$, and the difference
$\pile(a_{i+1}b_{i+1})∖\pile(a_ib_i)$ is $\{a_{i+1}▁\}$, an element of
the left column of $P$; 2) the step from $\aibi$ to $\aiibii$ is a
movement {\sl northeast} in the ZHA, as in from 34 to 35; here
$a_{i+1}b_{i+1} = a_i(b_i+1)$, and the difference
$\pile(a_{i+1}b_{i+1})∖\pile(a_ib_i)$ is $\{▁b_{i+1}\}$, an element of
the right column of $P$.

The easiest way to see how to convert from a set of question marks to
its associated slashing and vice-versa is by looking at an example.
Let's take the structure $((P,A),Q) \; \squigbijbody \; (H,S)$ on
which we've been working and build a table that shows how each step of
the path $(a_0b_0, a_1b_1, \ldots a_{10}b_{10})$ is ``seen'' by the
set $Q$, by the equivalence relations $\eqQ$, $\eqS$, $\eqL$, $\eqR$,
and by the slashing $S$ written in short form. We get:
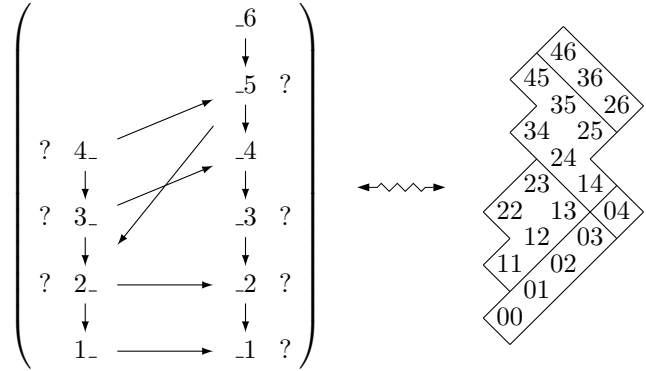
\def\Diffl#1#2#3#4#5#6#7#8{
  \pile(#3#4)∖\pile(#1#2)=\{#3▁\} &
  #3▁  #5\in Q &
  #1#2 #6\eqQ #3#4 &
  #1   #7\eqL #3   & &
  #3   #8     #1   & \\
  }
\def\Diffr#1#2#3#4#5#6#7#8{
  \pile(#3#4)∖\pile(#1#2)=\{▁#4\} &
  ▁#4  #5\in Q &
  #1#2 #6\eqQ #3#4 &
  & #2 #7\eqR   #4 &
  & #2 #8       #4 \\
  }
\def\Diffrsame #1#2#3#4{\Diffr #1#2 #3#4 {}   {}   {}     {}}
\def\Diffrother#1#2#3#4{\Diffr #1#2 #3#4 \not \not {\not} {∖}}
\def\Difflsame #1#2#3#4{\Diffl #1#2 #3#4 {}   {}   {}     {}}
\def\Difflother#1#2#3#4{\Diffl #1#2 #3#4 \not \not {\not} {/}}
$$\begin{array}{c}
  \tcg{(P,A),Q} \;\; \squigbij \;\;\; \zha{O_A(P),J}
  \\
  \\
  (a_0b_0, \ldots a_{10}b_{10}) = (00, 01, 02, 03, 04, 14, 24, 34, 35, 36, 46)
  \\
  \\
  \begin{array}{llllllll}
  \Difflsame  36 46 
  \Diffrother 35 36
  \Diffrsame  34 35
  \Difflsame  24 34
  \Difflsame  14 24
  \Difflother 04 14
  \Diffrother 03 04
  \Diffrsame  02 03
  \Diffrsame  01 02
  \Diffrsame  00 01
  \end{array}
  \end{array}
$$

There is an obvious correspondence between the elements of $P$ that
are {\sl not} in $Q$ and the `$/$'s and `$∖$' in $S$ that indicate
changes of equivalence class: $P∖Q = \{1▁, \; ▁4, ▁5\}$ corresponds to
$1/0$, $3∖4$, $5∖6$.

}

\directlua{tf_pop()}

% Local Variables:
% coding: utf-8-unix
% ee-tla: "jos"
% End:
    % (find-LATEX "2019J-ops-slashings.tex")
\newpage
% (find-LATEX "2019J-ops-logic.tex")
% (defun c () (interactive) (find-LATEXsh "lualatex -record 2019J-ops-logic.tex" :end))
% (defun d () (interactive) (find-pdf-page "~/LATEX/2019J-ops-logic.pdf"))
% (defun e () (interactive) (find-LATEX "2019J-ops-logic.tex"))
% (defun u () (interactive) (find-latex-upload-links "2019J-ops-logic"))
% (find-pdf-page   "~/LATEX/2019J-ops-logic.pdf")
% (find-sh0 "cp -v  ~/LATEX/2019J-ops-logic.pdf /tmp/")
% (find-sh0 "cp -v  ~/LATEX/2019J-ops-logic.pdf /tmp/pen/")
%   file:///home/edrx/LATEX/2019J-ops-logic.pdf
%               file:///tmp/2019J-ops-logic.pdf
%           file:///tmp/pen/2019J-ops-logic.pdf
% http://angg.twu.net/LATEX/2019J-ops-logic.pdf
% (find-LATEX "2019.mk")

\directlua{tf_push("2019J-ops-logic.tex")}

%      _                       
%     | |       ___  _ __  ___ 
%  _  | |_____ / _ \| '_ \/ __|
% | |_| |_____| (_) | |_) \__ \
%  \___/       \___/| .__/|___/
%                   |_|        
%
% «J-ops-and-regions» (to ".J-ops-and-regions")
% «J-operators»  (to ".J-operators")
% J-regions and J-operators
% (p2lp 1 "J-operators")
% (p2l    "J-operators")
\section{J-operators}
\label  {J-operators}
% Good (ph2p 9 "J-ops-and-regions")

% (fooi "Ω" "H")

A {\sl J-operator} on a Heyting Algebra $H ≡ (H,≤,⊤,⊥,∧,∨,→,↔,¬)$ is a
function $J:H→H$ that obeys the axioms $\J1$, $\J2$, $\J3$ below; we
usually write $J$ as $·^*:H→H$, and write the axioms as rules.
%
%L addabbrevs(".\\eqJ.", " \\eqJ ")
%L addabbrevs("&", "\\&", "vv", "∨", "->", "→")
%L -- addabbrevs("<=", "≤", "!!", "^{**}", "!", "^*")
% (fooi "!!" "²" "!" "¹" "^*" "¹" "<=" "≤" "->" "→" "&" "∧" "vv" "∨")
%:
%:    -----\J1   ------\J2   ------------\J3
%:    P≤P¹       P¹=P²       (P∧Q)¹=P¹∧Q¹
%:
%:    ^J1         ^J2          ^J3
%:
\pu
$$\ded{J1} \qquad \ded{J2} \qquad \ded{J3}$$
\par $\J1$ says that the operation $·^*$ is non-decreasing.
\par $\J2$ says that the operation $·^*$ is idempotent.
\par $\J3$ is a bit mysterious but will have interesting consequences.

\msk

% Note that when $H$ is a ZHA then any slash-operator on $H$ is a
% J-operator on it; see secs.\ref{slash-ops} and
% \ref{slash-ops-property}.
% 
% \msk

A J-operator induces an equivalence relation and equivalence classes
on $H$, like slashings do:
$$\begin{array}{rcl}
  P \eqJ Q  &\text{iff}& P^*=Q^* \\[5pt] \relax
  [P]^J &:=&         \setofst{Q∈H}{P^*=Q^*} \\
  \end{array}
$$
The equivalence classes of a J-operator $J$ are called {\sl
  $J$-regions}.

\msk

The axioms $\J1$, $\J2$, $\J3$ have many consequences. The first ones
are listed in Figure \ref{fig:J-ops-basic-derived-rules} as derived
rules, whose names mean:

$\Mop$ (monotonicity for products): a lemma used to prove $\Mo$,

$\Mo$ (monotonicity): $P≤Q$ implies $P^*≤Q^*$,

$\Sand$ (sandwiching): all truth values between $P$ and $P^*$ are equivalent,

$\ECa$: equivalence classes are closed by `$\&$',

$\ECv$: equivalence classes are closed by  `$∨$',

$\ECS$: equivalence classes are closed by sandwiching,

\begin{figure}
\pu
  \resizebox{\textwidth}{!}{%
  \fbox{$
  \def\bk{HELLO}
  \def\bk{\hspace{-1cm}}
  \begin{array}{rcl} \\
  \ded{Mop1}  &:=& \ded{Mop2}  \\ \\
  \ded{Mo1}   &:=& \ded{Mo2}   \\ \\
  \ded{Sand1} &:=& \ded{Sand2} \\ \\
  \ded{ECa1}  &:=& \ded{ECa2}  \\
  \ded{ECv1}  &:=& \ded{ECv2} \\ \\
  \ded{ECS1}  &:=& \ded{ECS2}  \\ \\
  \end{array}
  $}
  }
  \caption{J-operators: basic derived rules}
  \label{fig:J-ops-basic-derived-rules}
\end{figure}

\bsk

Take a J-equivalence class, $[P]^J$, and list its elements: $[P]^J =
\{P_1, \ldots, P_n\}$. Let $P_∧ := ((P_1∧P_2)∧\ldots)∧P_n$ and $P_∨ :=
((P_1∨P_2)∨\ldots)∨P_n$. Clearly $P_∧ ≤ P_i ≤ P_∨$ for each $i$, so
$[P]^J ⊆ [P_∧,P_∨]$. We will use the interval notation $[P,R]$ to mean
the set of all elements of $H$ obeying $P≤Q≤R$:
$$[P,R] \;\; = \;\; \setofst{Q∈H}{P≤Q≤R}.$$

Using $\ECa$ and $\ECv$ several times we see that:
$$\begin{array}{rrr}
                 P_1∧P_2 \eqJ P &&                P_1∨P_2 \eqJ P \\
           (P_1∧P_2)∧P_3 \eqJ P &&          (P_1∨P_2)∨P_3 \eqJ P \\
         \vdots\phantom{mmmm}   &&        \vdots\phantom{mmmm}   \\
  ((P_1∧P_2)∧\ldots)∧P_n \eqJ P && ((P_1∨P_2)∨\ldots)∨P_n \eqJ P \\
                     P_∧ \eqJ P &&                    P_∨ \eqJ P \\[5pt]
                    P_∧ ∈ [P]^J &&                   P_∨ ∈ [P]^J \\
  \end{array}
$$
and using $\ECS$ we can see that all elements between $P_∧$ and $P_∨$
are $J$-equivalent to $P$:
%:     
%:               P_∧.\eqJ.P   P_∨.\eqJ.P
%:               ----------   ----------
%:               {P_∧}^*=P^*  {P_∨}^*=P^*  
%:                ----------------------  
%:  P_∧≤Q≤P_∨        {P_∧}^*={P_∨}^*        
%:  --------------------------------\ECS
%:        {P_∧}^*=Q^*={P_∨}^*             {P_∨}^*=P^*
%:        -------------------------------------------
%:                  Q^*=P^*
%:                 ---------
%:                 Q.\eqJ.P
%:
%:                  ^foo
%:
$$\pu
  \ded{foo}
$$
so $[P_∧,P_∨] ⊆ [P]^J$. This means that {\sl J-regions are intervals}.

% (find-LATEX "2015planar-has.tex" "J-operators")
% (find-planarhaspage 13 "Part 2:" "J-operators and ZQuotients")
% (find-planarhastext 13 "Part 2:" "J-operators and ZQuotients")
% (find-LATEX "2015planar-has.tex" "J-derived-rules")
% (find-planarhaspage 15 "Derived rules")
% (find-planarhastext 15 "Derived rules")

\directlua{tf_pop()}

% Local Variables:
% coding: utf-8-unix
% ee-tla: "jol"
% End:
        % (find-LATEX "2019J-ops-logic.tex")
\newpage
% (find-LATEX "2019J-ops-midway.tex")
% (defun c () (interactive) (find-LATEXsh "lualatex -record 2019J-ops-midway.tex" :end))
% (defun d () (interactive) (find-pdf-page "~/LATEX/2019J-ops-midway.pdf"))
% (defun e () (interactive) (find-LATEX "2019J-ops-midway.tex"))
% (defun u () (interactive) (find-latex-upload-links "2019J-ops-midway"))
% (find-pdf-page   "~/LATEX/2019J-ops-midway.pdf")
% (find-sh0 "cp -v  ~/LATEX/2019J-ops-midway.pdf /tmp/")
% (find-sh0 "cp -v  ~/LATEX/2019J-ops-midway.pdf /tmp/pen/")
%   file:///home/edrx/LATEX/2019J-ops-midway.pdf
%               file:///tmp/2019J-ops-midway.pdf
%           file:///tmp/pen/2019J-ops-midway.pdf
% http://angg.twu.net/LATEX/2019J-ops-midway.pdf
% (find-LATEX "2019.mk")

% «.cuts-stopping-midway»		(to "cuts-stopping-midway")
% «.no-Y-cuts-no-L-cuts»		(to "no-Y-cuts-no-L-cuts")

\directlua{tf_push("2019J-ops-midway.tex")}

%   ____      _             _                    _             
%  / ___|   _| |_ ___   ___| |_ ___  _ __  _ __ (_)_ __   __ _ 
% | |  | | | | __/ __| / __| __/ _ \| '_ \| '_ \| | '_ \ / _` |
% | |__| |_| | |_\__ \ \__ \ || (_) | |_) | |_) | | | | | (_| |
%  \____\__,_|\__|___/ |___/\__\___/| .__/| .__/|_|_| |_|\__, |
%                                   |_|   |_|            |___/ 
%
% «cuts-stopping-midway»  (to ".cuts-stopping-midway")
% (p2lp 3 "cuts-stopping-midway")
% (p2l    "cuts-stopping-midway")
\section{Cuts stopping midway}
\label{cuts-stopping-midway}

Look at the figure below, that shows a partition of a ZHA $A=[00,66]$
into five regions, each region being an interval; this partition does
not come from a slashing, as it has cuts that stop midway. Define an
operation `$·^*$' on $A$, that works by taking each truth-value $P$ in
it to the top element of its region; for example, $30^*=61$.
%
%L mp = mpnew({def="foo", meta="10pt"}, "1234567654321")
%L mp:addlrs():addcuts("c 10w-14n 11n-61n 42w-46n 44n-04e"):output()
$$\pu
  \foo
$$
It is easy to see that `$·^*$' obeys $\J1$ and $\J2$; however, it does
{\sl not} obey $\J3$ --- we will prove that in
sec.\ref{no-Y-cuts-no-L-cuts}. As we will see, {\sl the partitions of
a ZHA into intervals that obey $\J1$, $\J2$, $\J3$ ae exactly the
slashings;} or, in other words, {\sl every J-operator comes from a
slashing}.

%  _   _        __   __             _       
% | \ | | ___   \ \ / /   ___ _   _| |_ ___ 
% |  \| |/ _ \   \ V /   / __| | | | __/ __|
% | |\  | (_) |   | |   | (__| |_| | |_\__ \
% |_| \_|\___/    |_|    \___|\__,_|\__|___/
%                                           
% «no-Y-cuts-no-L-cuts» (to ".no-Y-cuts-no-L-cuts")
% (p2lp 4 "no-Y-cuts-no-L-cuts")
% (p2l    "no-Y-cuts-no-L-cuts")
\subsection{The are no Y-cuts and no $\lambda$-cuts}
\label  {no-Y-cuts-no-L-cuts}
% Good (ph2p 12 "no-Y-cuts-no-L-cuts")

We want to see that if a partition of a ZHA $H$ into intervals has
``Y-cuts'' or ``$λ$-cuts'', like these parts of the last diagram in
sec.\ref{cuts-stopping-midway},
%
%R local Y, La =
%R 2/  22  \, 2/  25  \
%R  |21  12|   |24  15|
%R  \  11  /   \  14  /
%R 
%R Y :tozmp({def="Ycut", scale="10pt"}):addcells():addcuts("00w-01n 10e-10n"):output()
%R La:tozmp({def="Lcut", scale="10pt"}):addcells():addcuts("00w-00n 00e-10n"):output()
$$\pu
  \begin{array}{rcl}
  \Ycut &\Leftarrow& \text{this is a Y-cut}  \\\\
  \Lcut &\Leftarrow& \text{this is a $λ$-cut}\\
  \end{array}
$$
then the operation $J$ that takes each element to the top of its
equivalence class cannot obey $\J1$, $\J2$ and $\J3$ at the same time.
We will prove that by deriving rules that say that if $11 \eqJ 12$
then $21 \eqJ 22$, and that if $15 \eqJ 25$ then $14 \eqJ 24$;
actually, our rules will say that if $11^* = 12^*$ then $(11∨21)^* =
(12∨21)^*$, and that if $15^*=25^*$ then $(15∧24)^*=(25∧24)^*$. The
rules are:

% (fooi "!!" "²" "!" "¹" "^*" "¹" "<=" "≤" "->" "→" "&" "∧" "vv" "∨")
%
%:                                     Q¹=R¹
%:                                  ---------
%:                                  P∨Q¹=P∨R¹
%:                               ---------------
%:      Q¹=R¹                    (P∨Q¹)¹=(P∨R¹)¹
%:  ---------------\NoYcuts  :=  ===============\ostarcube
%:  (P∨Q)¹=(P∨R)¹                 (P∨Q)¹=(P∨R)¹           
%:
%:  ^NoYcuts1-                     ^NoYcuts2-
%:
%:                                     P¹=Q¹
%:                                  ---------
%:                                  P∨R¹=Q∨R¹
%:                               ---------------
%:      P¹=Q¹                    (P∨R¹)¹=(Q∨R¹)¹
%:  ---------------\NoYcuts  :=  ================\oor_6=\oor_4
%:  (P∨R)¹=(Q∨R)¹               (P∨R)¹=(Q∨R)¹
%:
%:  ^NoYcuts1                      ^NoYcuts2
%:
%:
%:                                  Q¹=R¹
%:                               ----------
%:      Q¹=R¹                    P¹∧Q¹=P¹∧R¹
%:  ---------------\NoLcuts  :=  ------------\J3
%:  (P∧Q)¹=(P∧R)¹                (P∧Q)¹=(P∧R)¹
%:
%:  ^NoLcuts1-                    ^NoLcuts2-
%:
%:                                  P¹=Q¹
%:                               ----------
%:      P¹=Q¹                    P¹∧R¹=Q¹∧R¹
%:  ---------------\NoLcuts  :=  ------------\J3
%:  (P∧R)¹=(Q∧R)¹                (P∧R)¹=(Q∧R)¹
%:
%:  ^NoLcuts1                     ^NoLcuts2
%:
$$\pu
  \begin{array}{rcl}
  \ded{NoYcuts1} &:=& \ded{NoYcuts2} \\ \\
  \ded{NoLcuts1} &:=& \ded{NoLcuts2} \\ \\
  \end{array}
$$

The expansion of double bar labeled `$\oor_6=\oor_4$' in the top
derivation uses twice the derived rule $\oor_6=\oor_4$, that is easy
to prove using the cubes of sec.\ref{cubes}.

\directlua{tf_pop()}

% Local Variables:
% coding: utf-8-unix
% ee-tla: "jom"
% End:
       % (find-LATEX "2019J-ops-midway.tex")
\newpage
% (find-LATEX "2019J-ops-cubes.tex")
% (defun c () (interactive) (find-LATEXsh "lualatex -record 2019J-ops-cubes.tex" :end))
% (defun d () (interactive) (find-pdf-page "~/LATEX/2019J-ops-cubes.pdf"))
% (defun e () (interactive) (find-LATEX "2019J-ops-cubes.tex"))
% (defun u () (interactive) (find-latex-upload-links "2019J-ops-cubes"))
% (find-pdf-page   "~/LATEX/2019J-ops-cubes.pdf")
% (find-sh0 "cp -v  ~/LATEX/2019J-ops-cubes.pdf /tmp/")
% (find-sh0 "cp -v  ~/LATEX/2019J-ops-cubes.pdf /tmp/pen/")
%   file:///home/edrx/LATEX/2019J-ops-cubes.pdf
%               file:///tmp/2019J-ops-cubes.pdf
%           file:///tmp/pen/2019J-ops-cubes.pdf
% http://angg.twu.net/LATEX/2019J-ops-cubes.pdf
% (find-LATEX "2019.mk")

% «.cubes»				(to "cubes")

\directlua{tf_push("2019J-ops-cubes.tex")}

%   ___  _          _                              _               
%  / _ \| |____   _(_) ___  _   _ ___    ___ _   _| |__   ___  ___ 
% | | | | '_ \ \ / / |/ _ \| | | / __|  / __| | | | '_ \ / _ \/ __|
% | |_| | |_) \ V /| | (_) | |_| \__ \ | (__| |_| | |_) |  __/\__ \
%  \___/|_.__/ \_/ |_|\___/ \__,_|___/  \___|\__,_|_.__/ \___||___/
%                                                                  
% «cubes»  (to ".cubes")
% (ph2p 13 "obvious-cubes")
% (p2lp 4  "cubes")
% (p2l     "cubes")
\section{How J-operators interact with connectives}
\label  {cubes}

% (find-LATEX "2015planar-has.tex" "J-connectives")
% (find-planarhaspage 16 "How J-operators interact with the connectives")
% (find-planarhastext 16 "How J-operators interact with the connectives")

The axiom $\J3$ says that $(P∧Q)¹=P¹∧Q¹$ --- it says something about
how `$·^*$' interacts with `$∧$'. Let's introduce a shorter notation.
There are eight ways to replace each of the `?'s in $(P^? ∧ Q^?)^?$ by
either nothing or a star. We establish that the three `?'s in $(P^? ∧
Q^?)^?$ are ``worth'' 1, 2 and 4 respectively, and we use $P \oand_n
Q$ to denote $(P^? ∧ Q^?)^?$ with the bits ``that belong to $n$''
replaced by stars. So:
$$\begin{array}{rcrcrcr}
  \oand_0     &=&   P∧Q,   && \oand_4 &=& (P∧Q)^*, \\
  \oand_1     &=& P^*∧Q,   && \oand_5 &=& (P^*∧Q)^*, \\
  \oand_2     &=&   P∧Q^*, && \oand_6 &=& (P∧Q^*)^*, \\
  \oand_3     &=& P^*∧Q^*, && \oand_7 &=& (P^*∧Q^*)^*. \\
  \end{array}
$$

We omit the arguments of $\oand_n$ when they are $P$ and $Q$ --- so we
can rewrite $(P∧Q)¹=P¹∧Q¹$ as $\oand_4=\oand_3$. These conventions
also hold for $\oor$ and $\oimp$.

\pu
%$$
%  \begin{array}{rcl}
%  \diag{cube-and*-obvious} && \diag{cube-and*-full} \\ \\
%  \diag{cube-or*-obvious}  && \diag{cube-or*-full}  \\ \\
%  \diag{cube-imp*-obvious} && \diag{cube-imp*-full} \\
%  \end{array}
%$$

It is easy to prove each one of the arrows in the cubes below ($A
\diagxyto/->/ B$ means $A≤B$):
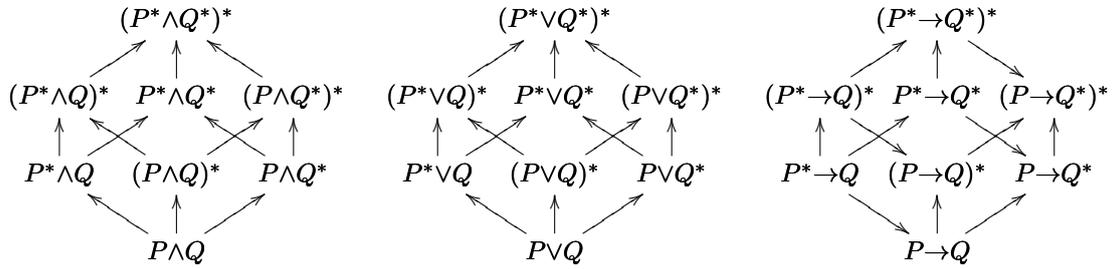
$$\myresizebox{$
  \pu
  \diag{cube-and*-obvious}
  \quad
  \diag{cube-or*-obvious}
  \quad
  \diag{cube-imp*-obvious}
$}
$$

\bsk

Let's write their sets of elements as $\oand_\zs := \{\oand_0, \ldots,
\oand_7\}$, $\oor_\zs := \{\oor_0, \ldots, \oor_7\}$, and $\oimp_\zs
:= \{\oimp_0, \ldots, \oimp_7\}$. The cubes above --- we will call
them the ``obvious and-cube'', the ``obvious or-cube'', and the
``obvious implication-cube'' --- can be interpreted as directed graphs
$(\oand_\zs, \OCube_\land)$, $(\oor_\zs, \OCube_\lor)$, $(\oimp_\zs,
\OCube_\to)$.

The ``extended cubes'' will be the directed graphs with the arrows
above plus the ones coming from these derived rules:
%:
%:
%:   =====================\oand_7=\oand_3=\oand_4
%:   (P¹∧Q¹)¹=P¹∧Q¹=(P∧Q)¹       
%:
%:   ^and*-extra-arrow           
%:
%:   ===============\oor_7≤\oor_3 
%:   (P¹∨Q¹)¹≤(P∨Q)¹       
%:                         
%:   ^or*-extra-arrow      
%:
%:   =============\oimp_6≤\oimp_3
%:   (P→Q¹)¹≤P¹→Q¹    
%:                    
%:   ^imp*-extra-arrow
%:
%:
%:
%:      ------\J2   ------\J2
%:      P²=P¹      Q²=Q¹
%:   =============================\J3
%:   (P¹∧Q¹)¹=P²∧Q²=P¹∧Q¹=(P∧Q)¹
%:   -----------------------------
%:      (P¹∧Q¹)¹=P¹∧Q¹=(P∧Q)¹
%:
%:      ^and*-extra-arrow-proof
%:
%:                                  ---------
%:                                  P→Q¹≤P→Q¹
%:      -----        -----        -----------
%:      P≤P∨Q        Q≤P∨Q        (P→Q¹)∧P≤Q¹
%:   ---------\Mo  ----------\Mo  ---------------\Mo
%:   P¹≤(P∨Q)¹     Q¹≤(P∨Q)¹       ((P→Q¹)∧P)¹≤Q²
%:   -----------------------       ---------------\J2
%:        P¹∨Q¹≤(P∨Q)¹              ((P→Q¹)∧P)¹≤Q¹
%:     ---------------\Mo          ---------------\J3
%:     (P¹∨Q¹)¹≤(P∨Q)²             (P→Q¹)¹∧P¹≤Q¹
%:     ----------------\J2         --------------
%:      (P¹∨Q¹)¹≤(P∨Q)¹             (P→Q¹)¹≤P¹→Q¹
%:
%:      ^or*-extra-arrow-proof         ^imp*-extra-arrow-proof
%:
%:
$$
\myresizebox{$\pu
  \def\bk{HELLO}
  \def\bk{\hspace{-0.5cm}}
  \begin{array}{rcl}
  \multicolumn{3}{l}{\ded{and*-extra-arrow} \quad :=} \\ \\
   \multicolumn{3}{r}{\ded{and*-extra-arrow-proof}}   \\ \\
  \ded{or*-extra-arrow}  &:=& \bk \ded{or*-extra-arrow-proof}  \\\\
  \ded{imp*-extra-arrow} &:=&     \ded{imp*-extra-arrow-proof} \\
  \end{array}
$}
$$
where $\oand_7=\oand_3=\oand_4$ will be interpreted as these arrows:
$$(P^*∧Q^*)^* \two/<-`->/<200> P^*∧Q^* \two/<-`->/<200> (P∧Q)^*$$

The directed graphs of these ``extended cubes'' will be called
$(\oand_\zs, \ECube_\land)$, $(\oor_\zs, \ECube_\lor)$, $(\oimp_\zs,
\ECube_\to)$. We are interested in the (non-strict) partial orders
that they generate, and we want an easy way to remember these partial
orders. The figure below shows these extended cubes at the left, and
at the right the ``simplified cubes'', $\SCube_∧$, $\SCube_∨$, and
$\SCube_→$, that generate the same partial orders that the extended
cubes.

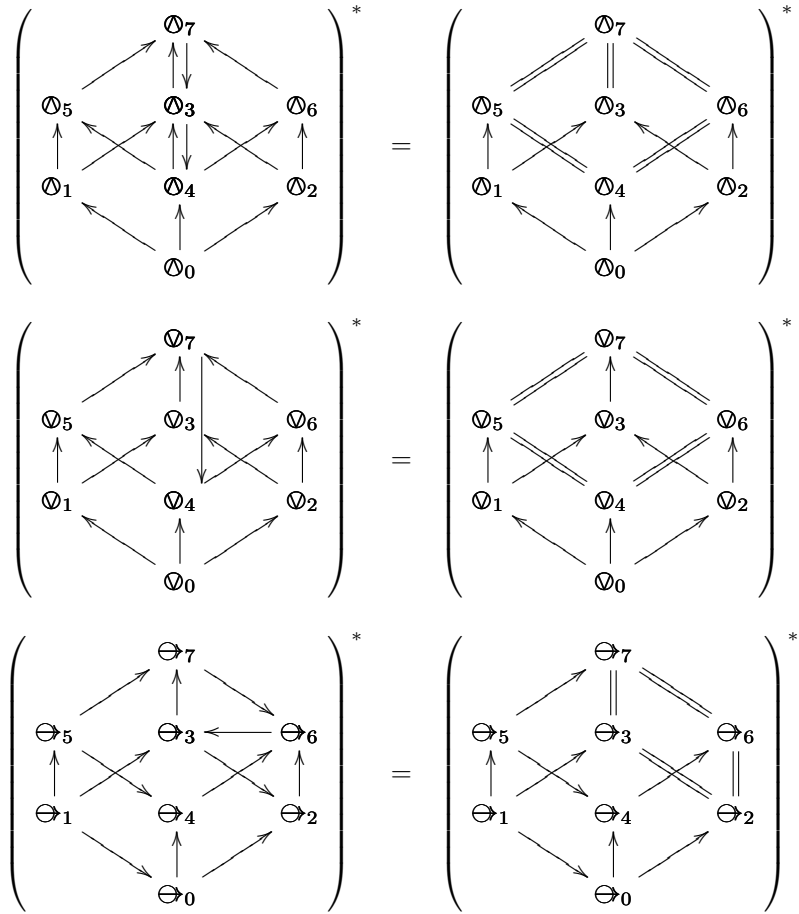
\pu
$$\resizebox{!}{180pt}{%
  $
    \begin{array}{rcl}
    \left( \diag{extended-and-cube-with-arrows} \right)^* &=&
    \left( \diag{extended-and-cube-with-equals} \right)^*     \\ \\
    \left( \diag{extended-or-cube-with-arrows}  \right)^* &=&
    \left( \diag{extended-or-cube-with-equals}  \right)^*     \\ \\
    \left( \diag{extended-imp-cube-with-arrows} \right)^* &=&
    \left( \diag{extended-imp-cube-with-equals} \right)^*     \\ \\
    \end{array}
  $
  }
$$
  %$}
%  \caption{The extended cubes and the simplified cubes.}
%  \label{fig:extended-and-simplified-cubes}
%\end{figure}

From these cubes it is easy to see, for example, that we can prove
$\oor_5 = \oor_6$ (as a derived rule).

\directlua{tf_pop()}

% Local Variables:
% coding: utf-8-unix
% ee-tla: "joc"
% End:
        % (find-LATEX "2019J-ops-cubes.tex")
\newpage
% (find-LATEX "2019J-ops-valuations.tex")
% (defun c () (interactive) (find-LATEXsh "lualatex -record 2019J-ops-valuations.tex" :end))
% (defun d () (interactive) (find-pdf-page "~/LATEX/2019J-ops-valuations.pdf"))
% (defun e () (interactive) (find-LATEX "2019J-ops-valuations.tex"))
% (defun u () (interactive) (find-latex-upload-links "2019J-ops-valuations"))
% (find-pdf-page   "~/LATEX/2019J-ops-valuations.pdf")
% (find-sh0 "cp -v  ~/LATEX/2019J-ops-valuations.pdf /tmp/")
% (find-sh0 "cp -v  ~/LATEX/2019J-ops-valuations.pdf /tmp/pen/")
%   file:///home/edrx/LATEX/2019J-ops-valuations.pdf
%               file:///tmp/2019J-ops-valuations.pdf
%           file:///tmp/pen/2019J-ops-valuations.pdf
% http://angg.twu.net/LATEX/2019J-ops-valuations.pdf
% (find-LATEX "2019.mk")

% «.valuations»				(to "valuations")

\directlua{tf_push("2019J-ops-valuations.tex")}

% __     __    _             _   _                 
% \ \   / /_ _| |_   _  __ _| |_(_) ___  _ __  ___ 
%  \ \ / / _` | | | | |/ _` | __| |/ _ \| '_ \/ __|
%   \ V / (_| | | |_| | (_| | |_| | (_) | | | \__ \
%    \_/ \__,_|_|\__,_|\__,_|\__|_|\___/|_| |_|___/
%                                                  
% «valuations»  (to ".valuations")
% (p2lp 7 "valuations")
% (p2l    "valuations")
\section{Valuations}
\label  {valuations}

Let $H_\odot$ and $J_\odot$ be a ZHA and a J-operator on it, and let
$v_\odot$ be a function from the set $\{P,Q\}$ to $H$. By an abuse of
language $v_\odot$ will also denote the triple $(H_\odot, J_\odot,
v_\odot)$ --- and by a second abuse of language $v_\odot$ will also
denote the obvious extension of $v_\odot: \{P,Q\}→H$ to the set of all
valid expressions formed from $P$, $Q$, $·^*$, $⊤$, $⊥$, and the
connectives.

Let $i,j∈\{0,\ldots,7\}$. Then $(\oand_i,\oand_j)∈\SCube^*_\land$
means that $\oand_i ≤ \oand_j$ is a theorem, and so $v_\odot(\oand_i)
≤ v_\odot(\oand_j)$ holds; i.e.,
$$\SCube^*_\land
  ⊆ \setofst {(\oand_i,\oand_j)}
             {i,j∈\{0,\ldots,7\}, \; v_\odot(\oand_i) ≤ v_\odot(\oand_j)}
$$
and the same for:
$$\begin{array}{c}
  \SCube^*_\lor
  ⊆ \setofst {(\oor_i,\oor_j)}
             {i,j∈\{0,\ldots,7\}, \; v_\odot(\oor_i) ≤ v_\odot(\oor_j)}
  \\
  \SCube^*_\to
  ⊆ \setofst {(\oimp_i,\oimp_j)}
             {i,j∈\{0,\ldots,7\}, \; v_\odot(\oimp_i) ≤ v_\odot(\oimp_j)}
  \\
  \end{array}
$$

Some valuations that turn these `$⊆$'s into `$=$'. Let
%
%L mp = mpnew({def="orCube", scale="11pt"}, "12321L"):addcuts("c 21/0 0|12")
%L mp:put(v"10", "P"):put(v"20", "P*", "P^*")
%L mp:put(v"01", "Q"):put(v"02", "Q*", "Q^*")
%L mp:output()
%
%L mp = mpnew({def="andCube", scale="11pt"}, "12321"):addcuts("c 2/10 01|2")
%L mp:put(v"20", "P"):put(v"21", "P*", "P^*")
%L mp:put(v"02", "Q"):put(v"12", "Q*", "Q^*")
%L mp:output()
%
%L mp = mpnew({def="impCube", scale="11pt"}, "12R1L"):addcuts("c 2/10 01|2")
%L mp:put(v"10", "P")  -- :put(v"20", "P*", "P^*")
%L mp:put(v"01", "Q")  -- :put(v"02", "Q*", "Q^*")
%L mp:output()
%
\pu
$$\begin{array}{c}
  (H_∧, J_∧, v_∧) = \andCube
  \qquad
  (H_∨, J_∨, v_∨) = \orCube \\
  (H_→, J_→, v_→) = \impCube \\
  \end{array}
$$
then
$$\begin{array}{c}
  \SCube^*_\land
  = \setofst {(\oand_i,\oand_j)}
             {i,j∈\{0,\ldots,7\}, \; v_∧(\oand_i) ≤ v_∧(\oand_j)}
  \\
  \SCube^*_\lor
  = \setofst {(\oor_i,\oor_j)}
             {i,j∈\{0,\ldots,7\}, \; v_∨(\oor_i) ≤ v_∨(\oor_j)}
  \\
  \SCube^*_\to
  = \setofst {(\oimp_i,\oimp_j)}
             {i,j∈\{0,\ldots,7\}, \; v_→(\oimp_i) ≤ v_→(\oimp_j)}
  \\
  \end{array}
$$
or, in more elementary terms:

\newpage

{\sl A very important fact.}
For any $i$ and $j$,
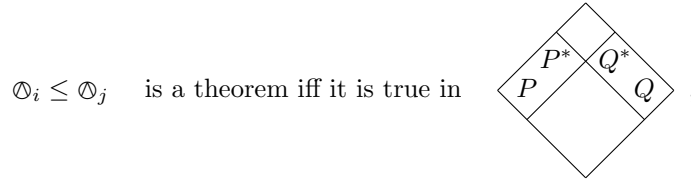
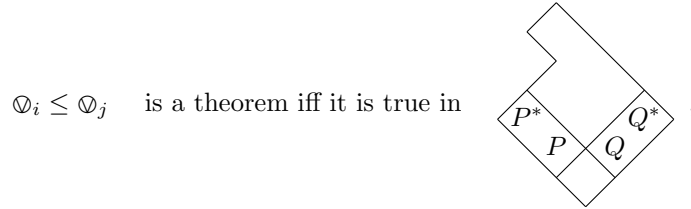
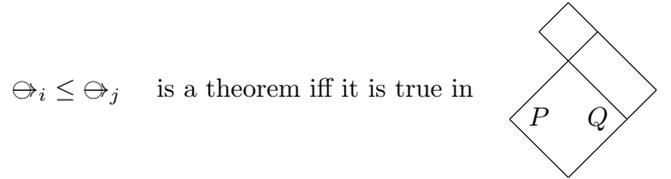
$$\pu
  \begin{array}{rcl}
  \oand_i≤\oand_j & \text{ is a theorem iff it is true in } & \andCube \;\; , \\
  \\
  \oor_i≤\oor_j & \text{ is a theorem iff it is true in } & \orCube  \;\; , \\
  \\
  \oimp_i≤\oimp_j & \text{ is a theorem iff it is true in } & \impCube \;\; . \\
  \end{array}
$$

The very important fact, and the valuations $v_∧$, $v_∨$, $v_→$, give
us:

\begin{itemize}

\item a way to {\sl remember} which sentences of the forms
  $\oand_i≤\oand_j$, $\oor_i≤\oor_j$, $\oimp_i≤\oimp_j$ are theorems;

\item countermodels for all the sentences of these forms not in
  $\SCube_∧$, $\SCube_∨$, $\SCube_→$. For example, $\oor_7≤\oor_4$ is
  not in $\SCube_∨$; and $v_∨(\oor_7)≤v_∨(\oor_4)$, which shows that
  $\oor_7≤\oor_4$ can't be a theorem.

\end{itemize}

% (find-books "__cats/__cats.el" "bell")
% (find-books "__cats/__cats.el" "bell" "163")

{\sl An observation.} I arrived at the cubes $\ECube_∧^*$,
$\ECube_∨^*$, $\ECube_→^*$ by taking the material in the corollary 5.3
of chapter 5 in \cite{BellLST} and trying to make it fit into less
mental space (as discussed in \cite{OchsIDARCT}); after that I wanted
to be sure that each arrow that is not in the extended cubes has a
countermodel, and I found the countermodels one by one; then I
wondered if I could find a single countermodel for all non-theorems in
$\ECube_∧^*$ (and the same for $\ECube_∨^*$ and $\ECube_→^*$), and I
tried to start with a valuation that distinguished {\sl some}
equivalence classes in $\ECube_∧^*$, and change it bit by bit, getting
valuations that distinguished more equivalence classes at every step.
Eventually I arrived at $v_∧$, $v_∨$ and at $v_→$, and at the ---
surprisingly nice --- ``very important fact'' above.

% (ph2p 20 "ZHA-vals-good")
% (ph2     "ZHA-vals-good")

Note that this valuation
%
%L mp = mpnew({def="orand", scale="11pt"}, "1234321L"):addcuts("c 432/10 01|23")
%L mp:put(v"20", "P"):put(v"31", "P*", "P^*")
%L mp:put(v"02", "Q"):put(v"13", "Q*", "Q^*")
%L mp:output()
%
$$(H_{∧∨},J_{∧∨},v_{∧∨}) \;\; = \;\; \pu\orand$$
distinguishes all equivalence classes in $\ECube^*_∧$ and in
$\ECube^*_∨$, but not in $\ECube^*_→$... it ``thinks'' that $P→Q$ and
$P^*→Q$ are equal.

\directlua{tf_pop()}

% Local Variables:
% coding: utf-8-unix
% ee-tla: "jov"
% End:
   % (find-LATEX "2019J-ops-valuations.tex")
\newpage
% (find-LATEX "2019J-ops-algebra.tex")
% (defun c () (interactive) (find-LATEXsh "lualatex -record 2019J-ops-algebra.tex" :end))
% (defun d () (interactive) (find-pdf-page "~/LATEX/2019J-ops-algebra.pdf"))
% (defun e () (interactive) (find-LATEX "2019J-ops-algebra.tex"))
% (defun u () (interactive) (find-latex-upload-links "2019J-ops-algebra"))
% (find-pdf-page   "~/LATEX/2019J-ops-algebra.pdf")
% (find-sh0 "cp -v  ~/LATEX/2019J-ops-algebra.pdf /tmp/")
% (find-sh0 "cp -v  ~/LATEX/2019J-ops-algebra.pdf /tmp/pen/")
%   file:///home/edrx/LATEX/2019J-ops-algebra.pdf
%               file:///tmp/2019J-ops-algebra.pdf
%           file:///tmp/pen/2019J-ops-algebra.pdf
% http://angg.twu.net/LATEX/2019J-ops-algebra.pdf
% (find-LATEX "2019.mk")

% «.polynomial-J-ops»			(to "polynomial-J-ops")
% «.algebra-of-piccs»			(to "algebra-of-piccs")
% «.algebra-of-J-ops»			(to "algebra-of-J-ops")
% «.slashings-are-poly»			(to "slashings-are-poly")

% (find-LATEX "2017planar-has-2.tex" "polynomial-J-ops")

\directlua{tf_push("2019J-ops-algebra.tex")}

%  ____       _             _                       
% |  _ \ ___ | |_   _      | |       ___  _ __  ___ 
% | |_) / _ \| | | | |  _  | |_____ / _ \| '_ \/ __|
% |  __/ (_) | | |_| | | |_| |_____| (_) | |_) \__ \
% |_|   \___/|_|\__, |  \___/       \___/| .__/|___/
%               |___/                    |_|        
%
% «polynomial-J-ops» (to ".polynomial-J-ops")
% (jopp 17 "polynomial-J-ops")
% (joa     "polynomial-J-ops")

\section{Polynomial J-operators}
\label  {polynomial-J-ops}
% (ph2p 29 "polynomial-J-ops")
% (ph2     "polynomial-J-ops")
% (phop 22)
% (pho "J-examples")
% (find-books "__cats/__cats.el" "fourman")
% (find-slnm0753page (+ 14 331)   "polynomial")

\def\Jnotnot{{(¬¬)}}
\def\JiiR   {{({→→}R)}}
\def\JoQ    {{(∨Q)}}
\def\JiR    {{({→}R)}}
\def\JfoQR  {{(∨Q∧{→}R)}}
\def\JmiQ   {({→→}Q∧{→}Q)}

{

It is not hard to check that for any Heyting Algebra $H$ and any
$Q,R∈H$ the operations $\Jnotnot$, $\ldots$, $\JfoQR$ below are
J-operators:
%
%$$\begin{array}{rclcrcl}
%      (¬¬) &:=& λP{:}H.¬¬P                  &&     (¬¬)(P) &=& ¬¬P \\
%    \JiiR  &:=& λP{:}H.((P{→}R){→}R)       &&   \JiiR(P) &=& (P{→}R){→}R \\
%      \JoQ &:=& λP{:}H.(P{∨}Q)              &&    \JoQ(P) &=& P∨Q \\
%      \JiR &:=& λP{:}H.(P{→}R)              &&    \JiR(P) &=& P{→}R\\
%    \JfoQR &:=& λP{:}H.((P{∨}Q) ∧ (P{→}R)) &&  \JfoQR(P) &=& (P{∨}Q)∧(P{→}R) \\
%  \end{array}
%$$
%
$$\begin{array}{rclcrcl}
      (¬¬)(P) &=& ¬¬P \\
     \JiiR(P) &=& (P{→}R){→}R \\
      \JoQ(P) &=& P∨Q \\
      \JiR(P) &=& P{→}R\\
    \JfoQR(P) &=& (P{∨}Q)∧(P{→}R) \\
  \end{array}
$$

Checking that they are J-operators means checking that each of them
obeys $\J1$, $\J2$, $\J3$. Let's define formally what are $\J1$, $\J2$
and $\J3$ ``for a given $F:H→H$'':
$$\begin{array}{rcc}
      \J1_F &:=&           (P ≤ F(P))     \\
      \J2_F &:=&        (F(P) = F(F(P))    \\
      \J3_F &:=&    (F(P∧P') = F(P)∧F(P')) \\
  \end{array}
$$
and:
$$\J123_F \quad := \quad \J1_F ∧ \J2_F ∧ \J3_F.$$

% (ph1p 18 "logic-in-HAs")
% (ph1     "logic-in-HAs")

Checking that $\Jnotnot$ obeys $\J1$, $\J2$, $\J3$ means proving
$\J123_\Jnotnot$ using only the rules from intuitionist logic from
section 10 of \cite{OchsPH1}; we will leave the proof of this, of and
$\J123_\JiiR$, $\J123_\JoQ$, and so on, to the reader.

\msk

The J-operator $\JfoQR$ is a particular case of building more complex
J-operators from simpler ones. If $J,K: H→H$, we define:
$$(J∧K) := λP{:}H.(J(P){∧}K(P))$$
it not hard to prove $\J123_{(J∧K)}$ from $\J123_J$ and $\J123_K$
using only the rules from intuitionistic logic.

\msk

The J-operators above are the first examples of J-operators in Fourman
and Scott's ``Sheaves and Logic'' (\cite{Fourman}); they appear in
pages 329--331, but with these names (our notation for them is at the
right):

(i) {\sl The closed quotient,}
$$J_a p = a ∨ p \qquad J_Q = \JoQ.$$

(ii) {\sl The open quotient,}
$$J^a p = a→p   \qquad J^R = \JiR.$$

(iii) {\sl The Boolean quotient}.
$$B_a p = (p→a)→a  \qquad B_R = \JiiR.$$

(iv) {\sl The forcing quotient}.
$$(J_a∧J^b)p = (a∨p)∧(b→p) \qquad (J_Q∧J^R) = \JfoQR.$$

(vi) {\sl A mixed quotient.}
$$(B_a∧J^a)p = (p→a)→p   \qquad (B_Q∧J^Q) = \JmiQ.$$

\msk

The last one is tricky. From the definition of $B_a$ and $J^a$ what we
have is
$$(B_a∧J^a)p = ((p→a)→a)∧(a→p),$$
but it is possible to prove
$$((p→a)→a)∧(a→p) \;\;↔\;\; ((p→a)→p)$$
intuitionistically.

The operators above are ``polynomials on $P,Q,R,→,∧,∨,⊥$'' in the
terminology of Fourman/Scott: ``If we take a polynomial in $→,∧,∨,⊥$,
say, $f(p,a,b,\ldots)$, it is a decidable question whether for all
$a,b,\ldots$ it defines a J-operator'' (p.331).

\msk

When I started studying sheaves I spent several years without any
visual intuition about the J-operators above. I was saved by ZHAs and
brute force --- and the brute force method also helps in testing if a
polynomial (in the sense above) is a J-operator in a particular case.
For example, take the operators $λP{:}H.(P∧22)$ and $({∨}22)$ on
$H=[00,44]$:
%
% (phop 23)
% (pho "J-ops-diagrams")
% (pho "J-ops-diagrams" "jout")
% (find-dn6 "zhas.lua" "shortoperators-tests")
% (find-dn6file "zhas.lua" "mpnewJ =")
%
%L shortoperators()
%L mpnewF = function (opts, spec, J)
%L     return mpnew(opts, spec, J):setz():zhaJ()
%L   end
%L
%L mpnewF({def="fooa"}, "123454321", function (P) return And(v"22", P) end):output()
%L mpnewF({def="fooo"}, "123454321", Cloq(v"22")):output()
%L mpnewJ({def="fooO"}, "123454321", Cloq(v"22")):zhaPs("22"):output()
%
$$\pu
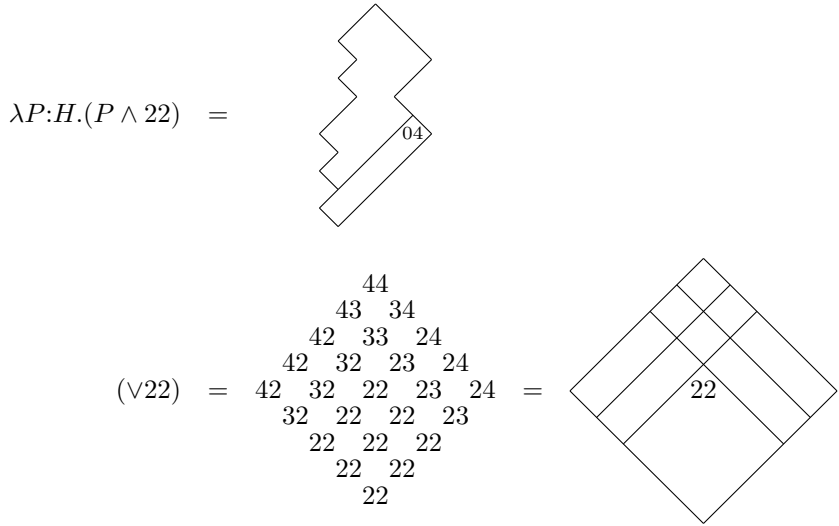
  \begin{array}{rcccl}
  λP{:}H.(P∧22) &=& \fooa \\ \\
        ({∨}22) &=& \fooo &=& \fooO \\
  \end{array}
$$

The first one, $λP{:}H.(P∧22)$, is not a J-operator; one easy way to
see that is to look at the region in which the result is 22 --- its
top element is 44, and this violates the conditions on slash-operators
in sec.\ref{slash-ops}. The second operator, $({∨}22)$, is a slash
operator and a J-operator; at the right we introduce a convenient
notation for visualizing the action of a polynomial slash-operator, in
which we draw only the contours of the equivalence classes and the
constants that appear in the polynomial.

Using this new notation, we have:
%
%L mpnewJ({def="fooboo", scale="7pt", meta="s"}, "123R2L1",   Booq(v"00")):zhaPs("00"):output()
%L mpnewJ({def="foobii", scale="7pt", meta="s"}, "123R2L1",   Booq(v"11")):zhaPs("11"):output()
%L
%L mpnewJ({def="fooboo", scale="7pt", meta="s"}, "123454321",   Booq(v"00")):zhaPs("00"):output()
%L mpnewJ({def="foobii", scale="7pt", meta="s"}, "123454321",   Booq(v"22")):zhaPs("22"):output()
%L mpnewJ({def="foobor", scale="7pt", meta="s"}, "1234567654321", Cloq(v"42")):zhaPs("42"):output()
%L mpnewJ({def="foobim", scale="7pt", meta="s"}, "1234567654321", Opnq(v"24")):zhaPs("24"):output()
%L mpnewJ({def="foofor", scale="7pt", meta="s"}, "1234567654321", Forq(v"42", v"24")):zhaPs("42 24"):output()
%L mpnewJ({def="foomix", scale="7pt", meta="s"}, "12345654321",   Mixq(v"22")):zhaPs("22"):output()
%
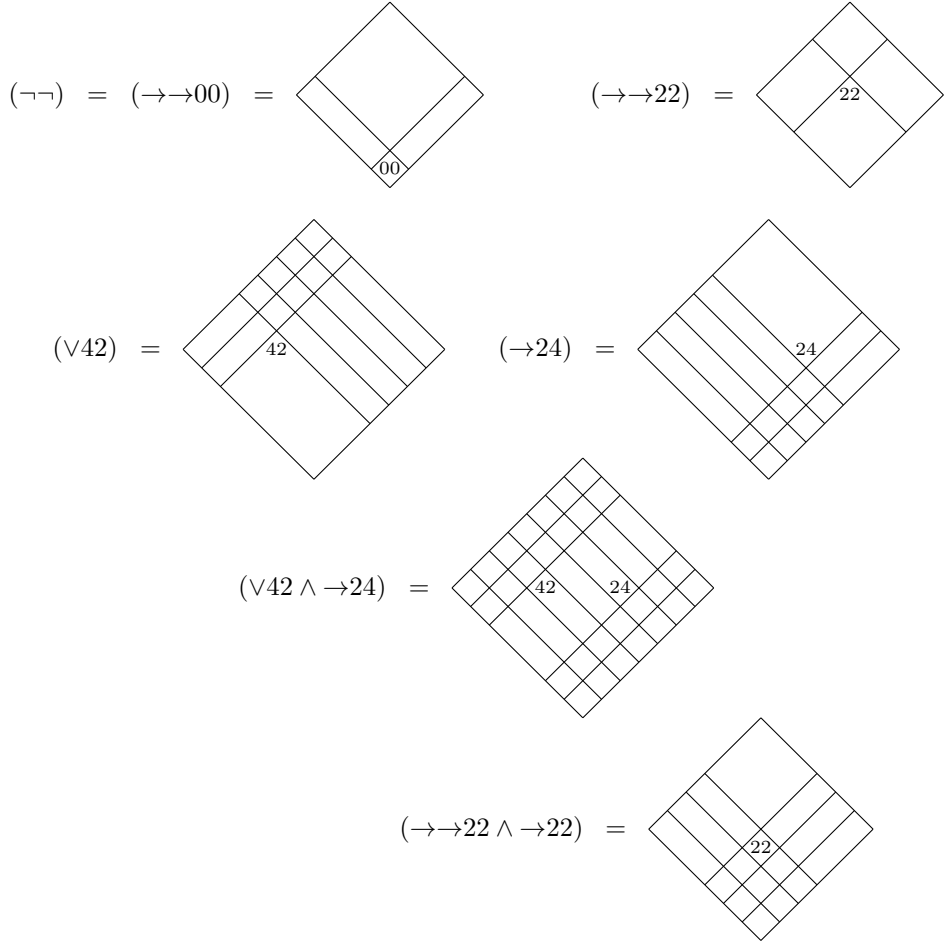
\pu
$$
  \begin{array}{c}
        (¬¬) \;\;=\;\;
        ({→→}00) \;\;=\;\; \fooboo \qquad \qquad
        ({→→}22) \;\;=\;\; \foobii \\
        \\
        ({∨}42) \;\;=\;\; \foobor \qquad
        ({→}24) \;\;=\;\; \foobim \\[-20pt]
        \\
        ({∨}42∧{→}24) \;\;=\;\; \foofor \\
        \qquad \qquad
        \qquad \qquad
        \qquad \qquad
        ({→→}22∧{→}22) \;\;=\;\; \foomix \\
  \end{array}
$$

Note that the slashing for $({∨}42∧{→}24)$ has all the cuts for
$({∨}42)$ plus all the cuts for $({→}24)$, and $({∨}42∧{→}24)$
``forces $42≤24$'' in the following sense: if $P^* = ({∨}42∧{→}24)(P)$
then $42^*≤24^*$.

}

%        _                      _            _               
%  _ __ (_) ___ ___ ___    __ _| | __ _  ___| |__  _ __ __ _ 
% | '_ \| |/ __/ __/ __|  / _` | |/ _` |/ _ \ '_ \| '__/ _` |
% | |_) | | (_| (__\__ \ | (_| | | (_| |  __/ |_) | | | (_| |
% | .__/|_|\___\___|___/  \__,_|_|\__, |\___|_.__/|_|  \__,_|
% |_|                             |___/                      
%
% «algebra-of-piccs» (to ".algebra-of-piccs")
\subsection{An algebra of piccs}
\label  {algebra-of-piccs}

We saw in the last section a case in which $(J∧K)$ has all the cuts
from $J$ plus all the cuts from $K$; this suggests that we {\sl may}
have an operation dual to that, that behaves as this: $(J∨K)$ has
exactly the cuts that are both in $J$ and in $K$:
$$\begin{array}{rcl}
  \Cuts(J∧K) &=& \Cuts(J)∪\Cuts(K) \\
  \Cuts(J∨K) &=& \Cuts(J)∩\Cuts(K) \\
  \end{array}
$$

And it $J_1, \ldots, J_n$ are all the slash-operators on a given ZHA,
then
$$\begin{array}{rclcl}
  \Cuts(J_1∧\ldots∧J_n) &=& \Cuts(J_1)∪\ldots∪\Cuts(J_k) &=& \text{(all cuts)} \\
  \Cuts(J_1∨\ldots∨J_n) &=& \Cuts(J_1)∩\ldots∩\Cuts(J_k) &=& \text{(no cuts)} \\
  \end{array}
$$
yield the minimal element and the maximal element, respectively, of an
algebra of slash-operators; note that the slash-operator with ``all
cuts'' is the identity map $λP{:H}.P$, and the slash-operator with
``no cuts'' is the one that takes all elements to $⊤$: $λP{:H}.⊤$.
This yields a lattice of slash-operators, in which the partial order
is $J≤K$ iff $\Cuts(J) ⊇ \Cuts(K)$. This is somewhat counterintuitive
if we think in terms of cuts --- the order seems to be reversed ---
but it makes a lot of sense if we think in terms of piccs
(sec.\ref{piccs-and-slashings}) instead.

\msk

Each picc $P$ on $\{0,\ldots,n\}$ has an associated function $·^P$
that takes each element to the top element of its equivalence class.
If we define $P≤P'$ to mean $∀a∈\{0,\ldots,n\}.\,a^P≤a^{P'}$, then we
have this:
%
% (pho     "algebra-of-piccs")
% (phop 25 "algebra-of-piccs")
%
%L partitiongraph = function (opts, spec, ylabel)
%L     local mp = MixedPicture.new(opts)
%L     for y=0,5 do mp:put(v(-1, y), y.."") end
%L     for x=0,5 do mp:put(v(x, -1), x.."") end
%L     for a=0,5 do local aP = spec:sub(a+1, a+1)+0; mp:put(v(a, aP), "*") end
%L     mp.lp:addlineorvector(v(0, 0), v(6, 0), "vector")
%L     mp.lp:addlineorvector(v(0, 0), v(0, 6.5), "vector")
%L     mp:put(v(7, 0), "a")
%L     mp:put(v(0, 7), "aP", ylabel)
%L     return mp
%L   end
%L pg = function (def, spec, ylabel)
%L     return partitiongraph({def=def, scale="5pt", meta="s"}, spec, ylabel)
%L   end
%L
%L pg("grapha", "012345", "a^P"     ):output()
%L pg("graphb", "113355", "a^{P'}"  ):output()
%L pg("graphc", "115555", "a^{P''}" ):output()
%L pg("graphd", "555555", "a^{P'''}"):output()
%L
%
%  a^P                  a^P                  a^P                  a^P
%    ^                    ^                    ^                    ^
%  5 |         *        5 |       * *        5 |   * * * *        5 * * * * * *
%  4 |       *          4 |                  4 |                  4 |
%  3 |     *        <=  3 |   * *        <=  3 |              <=  3 |
%  2 |   *              2 |                  2 |                  2 |
%  1 | *                1 * *                1 * *                1 |
%  0 *----------> a     0 +----------> a     0 +----------> a     0 +----------> a
%    0 1 2 3 4 5          0 1 2 3 4 5          0 1 2 3 4 5          0 1 2 3 4 5
%
%     0|1|2|3|4|5   <=      01|23|45     <=       01|2345            012345
%
$$\pu
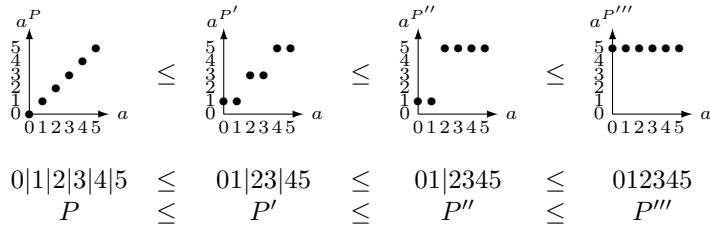
  \begin{array}{ccccccc}
  \grapha     &≤& \graphb  &≤& \graphc &≤& \graphd \\ \\
  0|1|2|3|4|5 &≤& 01|23|45 &≤& 01|2345 &≤& 012345  \\
  P           &≤& P'       &≤& P''     &≤& P'''    \\
  \end{array}
$$

This yields a partial order on piccs, whose bottom element is the
identity function $0|1|2|\ldots|n$, and the top element is $012\ldots
n$, that takes all elements to $n$.

The piccs on $\{0,\ldots,n\}$ form a Heyting Algebra, where
$⊤=01\ldots n$, $⊥=0|1|\ldots|n$, and `$∧$' and `$∨$' are the
operations that we have discussed above; it is possible to define a
`$→$' there, but this `$→$' is not going to be useful for us and we
are mentioning it just as a curiosity. We have, for example:
%
%D diagram algebra-of-piccs
%D 2Dx     100    +20   +20     +30 +20 +20
%D 2D  100       01234              T
%D 2D              ^                ^
%D 2D              |                |
%D 2D  +20      01|234             PvQ
%D 2D           ^    ^            ^   ^
%D 2D          /      \          /     \
%D 2D  +20 0|1|234  01|2|34     P       Q
%D 2D          ^      ^          ^     ^
%D 2D           \    /            \   /
%D 2D  +20     0|1|2|34            P&Q
%D 2D              ^                ^
%D 2D              |                |
%D 2D  +20     0|1|2|3|4           bot
%D 2D
%D (( T .tex= ⊤  PvQ .tex= P∨Q  P&Q .tex= P∧Q  bot .tex= ⊥
%D    01234 01|234 <-                             T PvQ <-
%D    01|234 0|1|234 <- 01|234 01|2|34 <-         PvQ P <- PvQ Q <-
%D    0|1|234 0|1|2|34 <- 01|2|34 0|1|2|34 <-     P P&Q <- Q P&Q <-
%D    0|1|2|34 0|1|2|3|4 <-                       P&Q bot <-
%D
%D ))
%D enddiagram
%D
$$\pu \diag{algebra-of-piccs}$$

%      _                               _            _               
%     | |       ___  _ __  ___    __ _| | __ _  ___| |__  _ __ __ _ 
%  _  | |_____ / _ \| '_ \/ __|  / _` | |/ _` |/ _ \ '_ \| '__/ _` |
% | |_| |_____| (_) | |_) \__ \ | (_| | | (_| |  __/ |_) | | | (_| |
%  \___/       \___/| .__/|___/  \__,_|_|\__, |\___|_.__/|_|  \__,_|
%                   |_|                  |___/                      
%
% «algebra-of-J-ops» (to ".algebra-of-J-ops")
\subsection{An algebra of J-operators}
\label  {algebra-of-J-ops}
% Bad (ph2p 50 "algebra-of-J-ops")

% (find-books "__cats/__cats.el" "fourman")

Fourman and Scott define the operations $∧$ and $∨$ on J-operators in
pages 325 and 329 (\cite{Fourman}), and in page 331 they list ten
properties of the algebra of J-operators:
$$
\def\li#1 #2 #3 #4    #5 #6 #7 #8 {\text{#1}&#2&#3&#4& &\text{#5}&#6&#7&#8& \\}
\def\li#1 #2 #3 #4                {\text{#1}&#2&#3&#4& }
\begin{array}{rlclcrlclc}
\li    (i) J_a∨J_b = J_{a∨b}  && (∨21)∨(∨12)=(∨22) \\
\li   (ii) J^a∨J^b = J^{a∧b}  && ({→}32)∨({→}23)=({→}22)  \\
\li  (iii) J_a∧J_b = J_{a∧b}  && (∨21)∧(∨12)=(∨11) \\
\li   (iv) J^a∧J^b = J^{a∨b}  && ({→}32)∧({→}23)=({→}33) \\
\li    (v) J_a∧J^a = ⊥        && (∨22)∧({→}22)=(⊥)  \\
\li   (vi) J_a∨J^a = ⊤        && (∨22)∨({→}22)=(⊤) \\
\li  (vii) J_a∨K   = K∘J_a      \\
\li (viii) J^a∨K   = J^a∘K    \\
\li   (ix) J_a∨B_a = B_a        \\
\li    (x) J^a∨B_b = B_{a→b}  \\
\end{array}
$$

% (pho     "J-ops-algebra-2")
% (phop 28 "J-ops-algebra-2")

The first six are easy to visualize; we won't treat the four last
ones. In the right column of the table above we've put a particular
case of (i), $\ldots$, (vi) in our notation, and the figures below put
all together.

In Fourman and Scott's notation,

%D diagram J-alg-1
%D 2Dx     100 +20 +20 +10 +10 +20 +20
%D 2D  100             T
%D 2D  +30     A               E
%D 2D  +20 B       C       F       G
%D 2D  +20     D               H
%D 2D  +30            bot
%D 2D
%D ren      T     ==>           J_⊤=⊤=J^⊥
%D ren    A   E   ==>     J_{22}          J^{22}
%D ren   B C F G  ==> J_{21} J_{12}   J^{32} J^{23}
%D ren    D   H   ==>     J_{11}          J^{11}
%D ren     bot    ==>           J_⊥=⊥=J^⊤
%D
%D (( A T -> E T ->
%D    B A -> C A -> F E -> G E ->
%D    D B -> D C -> H F -> H G ->
%D    bot D -> bot H ->
%D ))
%D enddiagram
%D
$$\pu
  \diag{J-alg-1}
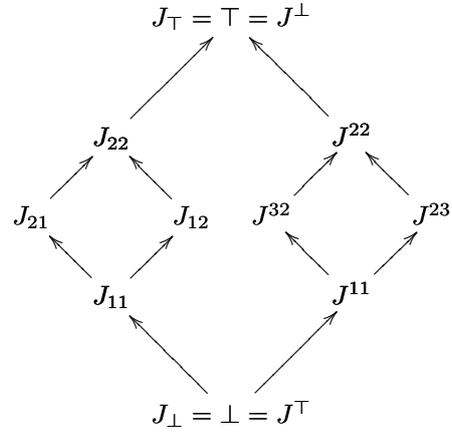
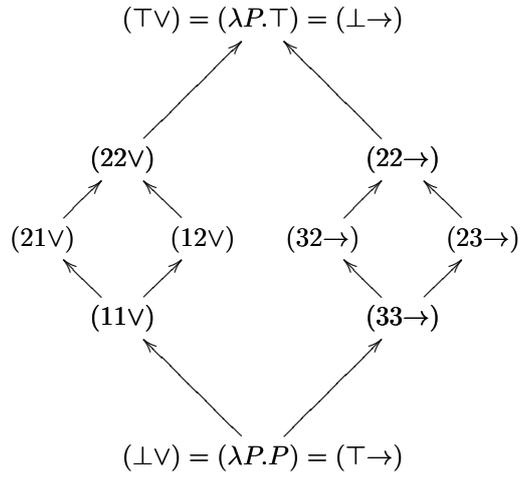
$$
in our notation,
%
%D diagram J-alg-2
%D 2Dx     100 +20 +20 +10 +10 +20 +20
%D 2D  100             T
%D 2D  +30     A               E
%D 2D  +20 B       C       F       G
%D 2D  +20     D               H
%D 2D  +30            bot
%D 2D
%D 2Dx     100 +20 +20 +15 +15 +20 +20
%D 2D  100             T
%D 2D  +35     A               E
%D 2D  +20 B       C       F       G
%D 2D  +20     D               H
%D 2D  +35            bot
%D 2D
%D ren      T     ==>      (⊤∨)=(λP.⊤)=(⊥{→})
%D ren    A   E   ==>     (22∨)          (22{→})
%D ren   B C F G  ==> (21∨) (12∨)   (32{→}) (23{→})
%D ren    D   H   ==>     (11∨)          (33{→})
%D ren     bot    ==>      (⊥∨)=(λP.P)=(⊤{→})
%D
%D (( A T -> E T ->
%D    B A -> C A -> F E -> G E ->
%D    D B -> D C -> H F -> H G ->
%D    bot D -> bot H ->
%D ))
%D enddiagram
%D
$$\pu
  \diag{J-alg-2}
$$
and drawing the polynomial J-operators as in
sec.\ref{polynomial-J-ops}:
%
%L deforp = function (P, draw, name)
%L     PP(P, name)
%L     mpnewJ({def=name, scale="4.5pt", meta="t"}, "123454321", Cloq(v(P))):zhaPs(draw):print():output()
%L   end
%L defimp = function (P, draw, name)
%L     PP(P, name)
%L     mpnewJ({def=name, scale="4.5pt", meta="t"}, "123454321", Opnq(v(P))):zhaPs(draw):print():output()
%L   end
%L deforp("21", "21", "ora")
%L deforp("22", "22", "orb")
%L deforp("11", "11", "orc")
%L deforp("12", "12", "ord")
%L deforp("00", "",   "orB")
%L defimp("32", "32", "ima")
%L defimp("22", "22", "imb")
%L defimp("33", "33", "imc")
%L defimp("34", "34", "imd")
%L defimp("00", "",   "imB")
%
%R local algebra =
%R 4/            !imB            \
%R  |                            |
%R  |    !orb            !imb    |
%R  |!ora    !ord    !ima    !imd|
%R  |    !orc            !imc    |
%R  |                            |
%R  \            !orB            /
%R 
%R algebra:tomp({def="foo", scale="30pt"}):addcells():output()
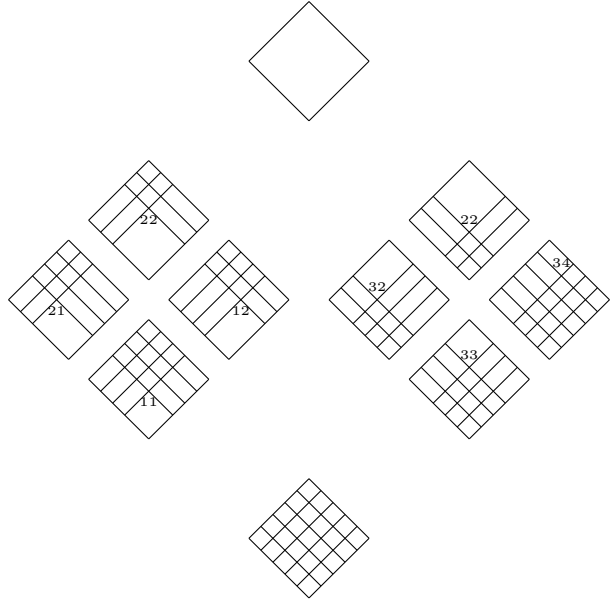
$$\pu
  % \bhbox{$
  \begin{array}{c}
  \\
  \foo \\
  \\
  \end{array}
  % $}
$$

% $$
% \def\li#1 #2 #3 #4    #5 #6 #7 #8 {\text{#1}&#2&#3&#4& &\text{#5}&#6&#7&#8& \\}
% \begin{array}{rlclcrlclc}
% \li   (i) J_a∨J_b = J_{a∨b}    (ii)   J^a∨J^b = J^{a∧b}
% \li (iii) J_a∧J_b = J_{a∧b}    (iv)   J^a∧J^b = J^{a∨b}
% \li   (v) J_a∧J^a = ⊥          (vi)   J_a∨J^a = ⊤
% \li (vii) J_a∨K   = K∘J_a      (viii) J^a∨K   = J^a∘K
% \li  (ix) J_a∨B_a = B_a        (x)    J^a∨B_b = B_{a→b}
% \end{array}
% $$
% 
% $$\begin{array}{rclcrcl}
%       (¬¬) &:=& λP:H.¬¬P                   &&     (¬¬)(P) &=& ¬¬P \\
%     (→→R) &:=& λP:H.((P{→}R){→}R)       &&   (→→R)(P) &=& (P{→}R){→}R \\
%       (∨Q) &:=& λP:H.(P{∨}Q)              &&    (∨Q)(P) &=& P∨Q \\
%       (→R) &:=& λP:H.(P{→}R)              &&     (→R)(P) &=& P{→}R\\
%   (∨Q∧→R) &:=& λP:H.((P{∨}Q) ∧ (P{→}R)) && (∨Q∧→R)(P) &=& (P{∨}Q)∧(P{→}R) \\
%   \end{array}
% $$

%  ____  _           _                   _       
% / ___|| | __ _ ___| |__    _ __   ___ | |_   _ 
% \___ \| |/ _` / __| '_ \  | '_ \ / _ \| | | | |
%  ___) | | (_| \__ \ | | | | |_) | (_) | | |_| |
% |____/|_|\__,_|___/_| |_| | .__/ \___/|_|\__, |
%                           |_|            |___/ 

% «slashings-are-poly» (to ".slashings-are-poly")
\subsection{All slash-operators are polynomial}
\label  {slashings-are-poly}
% (ph2p 51 "slashings-are-poly")

% (find-xpdfpage "~/LATEX/2015planar-has.pdf" 30)
% (find-LATEX "2015planar-has.tex" "zquotients-poly")

{

%L local ba, bb, bc = Booq(v"04"), Booq(v"23"), Booq(v"45")
%L local babc = Jand(ba, Jand(bb, bc))
%L mpnewJ({def="slaT", scale="7pt", meta="s"}, "1R2R3212RL1", babc       ):addlrs()         :output():print()
%L mpnewJ({def="slaA", scale="7pt", meta="s"}, "1R2R3212RL1", Booq(v"04")):addlrs()         :output():print()
%L mpnewJ({def="slaB", scale="7pt", meta="s"}, "1R2R3212RL1", Booq(v"23")):addlrs()         :output():print()
%L mpnewJ({def="slaC", scale="7pt", meta="s"}, "1R2R3212RL1", Booq(v"45")):addlrs()         :output():print()
%L mpnewJ({def="fooa", scale="7pt", meta="s"}, "1R2R3212RL1", Booq(v"04")):zhaPs("04")      :output():print()
%L mpnewJ({def="foob", scale="7pt", meta="s"}, "1R2R3212RL1", Booq(v"23")):zhaPs("23")      :output():print()
%L mpnewJ({def="fooc", scale="7pt", meta="s"}, "1R2R3212RL1", Booq(v"45")):zhaPs("45")      :output():print()
%L mpnewJ({def="food", scale="7pt", meta="s"}, "1R2R3212RL1", babc       ):zhaPs("04 23 45"):output():print()
%L mpnew ({def="fooe", scale="7pt", meta="s"}, "1R2R3212RL1"      ):addlrs():output()
%L mpnewJ({def="foof", scale="7pt", meta="s"}, "1R2R3212RL1", babc):zhaJ()  :output()
\pu

Here is an easy way to see that all slashings --- i.e., J-operators on
ZHAs --- are polynomial. Every slashing $J$ has only a finite number of
cuts; call them $J_1, \ldots, J_n$. For example:
$$J=\slaT \qquad J_1=\slaA \quad J_2=\slaB \quad J_3=\slaC$$

Each cut $J_i$ divides the ZHA into an upper region and a lower
region, and $J_i(00)$ yields the top element of the lower region.
Also, $({→→}J_i(00))$ is a polynomial way of expressing that cut:
$$\def\foo#1#2
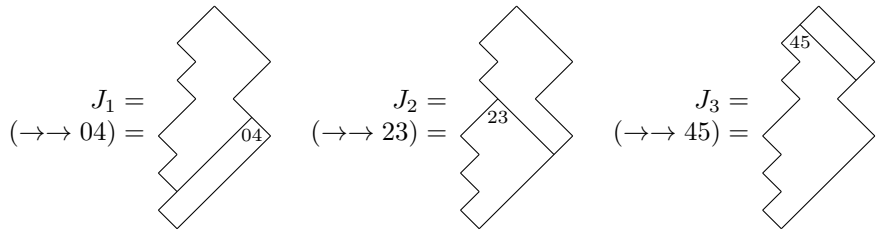
{
    \begin{array}{rr}
      J_#1= \\
      % ({→→}J_#1(00))= \\
      (→→#2)= \\
     \end{array}}
  \foo1{04} \fooa \quad
  \foo2{23} \foob \quad
  \foo3{45} \fooc
$$

The conjunction of these `$({→→}J_i(00))$'s yields the original
slashing:
%
% $$
%   \begin{array}{c}
%   \fooa ∧ \foob ∧ \fooc = \food
%   \end{array}
% $$
%
$$(→→04)∧(→→23)∧(→→45) = \food = J$$

}

\directlua{tf_pop()}

% Local Variables:
% coding: utf-8-unix
% ee-tla: "joa"
% End:
      % (find-LATEX "2019J-ops-algebra.tex")
\newpage
% (find-LATEX "2019J-ops-categories.tex")
% (defun c () (interactive) (find-LATEXsh "lualatex -record 2019J-ops-categories.tex" :end))
% (defun d () (interactive) (find-pdf-page "~/LATEX/2019J-ops-categories.pdf"))
% (defun e () (interactive) (find-LATEX "2019J-ops-categories.tex"))
% (defun u () (interactive) (find-latex-upload-links "2019J-ops-categories"))
% (find-pdf-page   "~/LATEX/2019J-ops-categories.pdf")
% (find-sh0 "cp -v  ~/LATEX/2019J-ops-categories.pdf /tmp/")
% (find-sh0 "cp -v  ~/LATEX/2019J-ops-categories.pdf /tmp/pen/")
%   file:///home/edrx/LATEX/2019J-ops-categories.pdf
%               file:///tmp/2019J-ops-categories.pdf
%           file:///tmp/pen/2019J-ops-categories.pdf
% http://angg.twu.net/LATEX/2019J-ops-categories.pdf
% (find-LATEX "2019.mk")

% «.classifier-big-figure-defs»	(to "classifier-big-figure-defs")
% «.classifier-intro»		(to "classifier-intro")
% «.classifier-big-figure»	(to "classifier-big-figure")
% «.Set-PA»			(to "Set-PA")
% «.Set-PA-logic»		(to "Set-PA-logic")
% «.Set-PA-morphisms»		(to "Set-PA-morphisms")

\directlua{tf_push("2019J-ops-categories.tex")}

\def\dnto{{\downarrow}}
\def\sfJ{\mathsf{J}}

%   ____      _       
%  / ___|__ _| |_ ___ 
% | |   / _` | __/ __|
% | |__| (_| | |_\__ \
%  \____\__,_|\__|___/
%                     
% «classifier-intro»  (to ".classifier-intro")
% (jopp 24 "classifier-intro")
% (joe     "classifier-intro")

\section{Categories, toposes, sheaves}

In this section I will explain {\sl very, very briefly} how to adapt
what we saw about J-operators to toposes. The first big diagram that
we will try to understand is the in one in
Figure \ref{fig:classifier-big} below, that shows in its upper part a
structure $((P,A),Q) \squigbij (H,J)$ with a 2CG with question marks
and its associated ZHA with J-operator, and in its lower part the
classifier $Ω$ of the topos $\Set^{(P,A)}$ and the local operator
$j:Ω→Ω$ that is associated to $Q$ and $J$. This big diagram shows how
to define the classifier and the local operator, and in
sec.\ref{kan-extensions} we will see another big diagram that shows
how to define sheaves and sheafification.

I will omit some technical details --- a very readable reference for
them is \cite{McLarty}, chapters 13 and 22. I learned most of them
from \cite{BellLST}, though.

\pu
\begin{figure}[h]
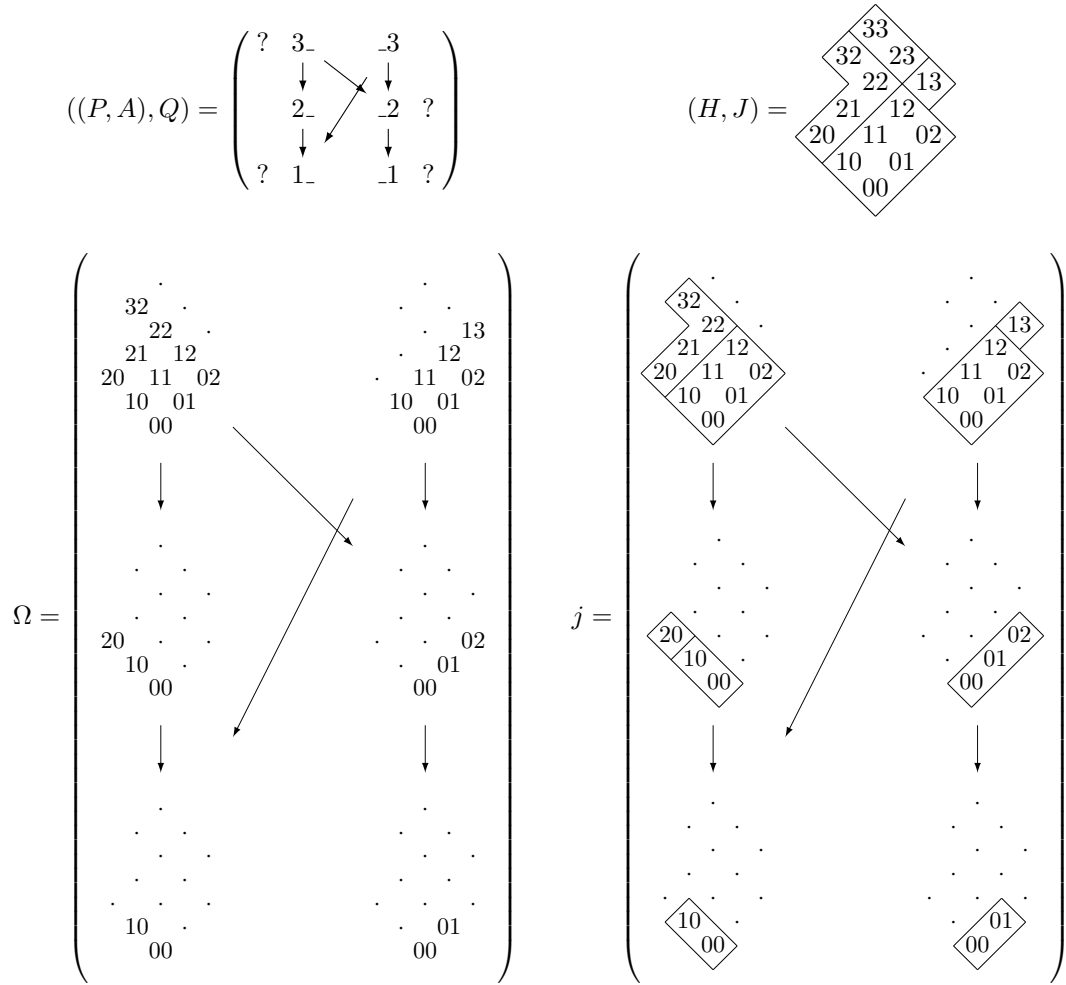

  $ % \hspace{-40pt}
  \begin{array}{ccc}
  ((P,A),Q) = \tcg{(P,A),Q 33}
    & 
    & (H,J) = \zha{H,J 33}
    \\ \\
  Ω = \scalebox{0.9}{$\tcg{Om}$} &
    &
    j = \scalebox{0.9}{$\tcg{j}$} \\
    \\
  \end{array}
  $
  \caption{The classifier and a local operator in a particular case}
  \label{fig:classifier-big}
\end{figure}

\newpage

%  ____       _     /\    ______     _  __  
% / ___|  ___| |_  |/\|  / /  _ \   / \ \ \ 
% \___ \ / _ \ __|      | || |_) | / _ \ | |
%  ___) |  __/ |_       | ||  __/ / ___ \| |
% |____/ \___|\__|      | ||_| ( )_/   \_\ |
%                        \_\   |/       /_/ 
%
% «Set-PA»  (to ".Set-PA")
% (jopp 23 "Set-PA")
% (joe     "Set-PA")

\subsection{Toposes of the form $\Set^{(P,A)}$}
\label{Set-PA}

%L kite  = ".1.|2.3|.4.|.5."
%L six   = ".1..|2.3.|.4.5|..6."
%L mp = MixedPicture.new({def="dagKite", meta="s", scale="5pt"}, z):zfunction(kite):output()
%L mp = MixedPicture.new({def="dagSix",  meta="s", scale="4.5pt"}, z):zfunction(six):output()
\pu

In sec.2 of \cite{OchsPH1} we established that the same bullet diagram
--- say, $\dagKite{•}{•}{•}{•}{•}$ --- could be intepreted as subset
of $\Z^2$, as a DAG, or as poset, depending on the context. Now we
will do something similar for graphs whose nodes are labeled. A
diagram like this
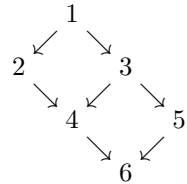
\pu
$$\zha{pB0}
$$
is interpreted as a DAG by default, but in this section it will be
also be interpreted as a (posetal) category in some contexts. We will
keep the same notation: if $(P,A)$ is a DAG then we will denote
$(P,A)$ ``regarded as a category'' by $(P,A)$.

A functor $F$ from a category $(P,A)$ to $\Set$ can be drawn as a
diagram with the same shape as $(P,A)$. If we draw the internal view
of $F:(P,A)→\Set$ over its internal view as in the introduction
of \cite{OchsPH1} we get this diagram:
%
%D diagram F6-int-ext
%D 2Dx     100    +75
%D 2D  100 I1 --> I2
%D 2D
%D 2D  +40 E1 --> E2
%D 2D
%D ren I1 I2 ==> \left(\zha{pB0}\right) \left(\zha{pF0}\right)
%D ren E1 E2 ==> (P,A)    \Set
%D
%D (( I1 I2 |->
%D    E1 E2  -> .plabel= a F
%D
%D ))
%D enddiagram
%D
$$\pu
  \diag{F6-int-ext}
$$
The `$\diagxyto/|->/$' in it stands for a bunch of
`$\diagxyto/|->/$'s, one for each object and one for each morphism.

We will only draw the upper-right part of diagrams like the one above.
With this convention, an object $F∈\Set^{(P,A)}$ can be drawn as:
$$
  F \;\; = \;\; \pBargs{F_1}{F_2}{F_3}{F_4}{F_5}{F_6}
$$

% (elep 6 "elephant-A2.1.3")
% (ele    "elephant-A2.1.3")

Every category of the form $\Set^{(P,A)}$ where $(P,A)$ is a finite
graph is a topos --- see \cite{EA}, example A2.1.3 --- so categories
of the form $\Set^{(P,A)}$ are toposes whose objects can be drawn as
$(P,A)$-shaped diagrams.

%  ____       _     /\    ______     _  __    _             _      
% / ___|  ___| |_  |/\|  / /  _ \   / \ \ \  | | ___   __ _(_) ___ 
% \___ \ / _ \ __|      | || |_) | / _ \ | | | |/ _ \ / _` | |/ __|
%  ___) |  __/ |_       | ||  __/ / ___ \| | | | (_) | (_| | | (__ 
% |____/ \___|\__|      | ||_| ( )_/   \_\ | |_|\___/ \__, |_|\___|
%                        \_\   |/       /_/           |___/        
%
% «Set-PA-logic»  (to ".Set-PA-logic")
% (jopp 26 "Set-PA-logic")
% (joe     "Set-PA-logic")

\subsection{The logic of toposes of the form $\Set^{(P,A)}$}

%L -- Define the components of the
%L -- internal views of Om and j.
%L
%L myspec = "1232RL1"
%L tdims_mini = TCGDims {h=32,  v=25,  q=15, crh=8,  crv=7, qrh=5}
%L tdims_big  = TCGDims {h=110, v=110,       crh=30, crv=45}
%L tspec_PAQ  = TCGSpec.new("33; 32, 13", "?.?", "??.")
%L tspec_PA   = TCGSpec.new("33; 32, 13")
%L tcg_PAQ    = TCGQ.newdsoa(tdims_mini, tspec_mini,
%L                           {tdef="(P,A),Q", meta="1pt p"},
%L                               "     h v p lr q"
%L                            -- "B QB h v p lr q"
%L                          )
%L tcg_j      = TCGQ.newdsoa(tdims_big, tspec_big,
%L                           {tdef="j", meta="1pt p"},
%L                           -- "B h v p"
%L                              "h v p"
%L                          )
%L tcg_Om     = TCGQ.newdsoa(tdims_big, tspec_big,
%L                           {tdef="Om", meta="1pt p"},
%L                           -- "B h v p"
%L                              "h v p"
%L                          )

The terminal object $1∈\Set^{(P,A)}$ is:
$$\def\u{\{*\}}
  1 \;\; = \;\; \pBargs{\u}{\u}{\u}{\u}{\u}{\u}
$$
and we can obtain all its subobjects by replacing some of the
`$\{*\}$'s in it by empty sets. If we rewrite each $\{*\}$ as 1 and
each $∅$ as 0 and use a more compact notation, then $1
= \dagSix111111$ and:
%
% (ph1p 25 "topologies-on-ZSets")
% (ph1     "topologies-on-ZSets")
%
$$\Sub(1) = \left\{
  \dagSix000000,
  \dagSix000001,
  \dagSix000011,
  \dagSix000101,
  \dagSix000111,
  \dagSix001111,
  \dagSix010101,
  \dagSix010111,
  \dagSix011111,
  \dagSix111111
  \right\}
$$
The Heyting Algebra of subobjects of 1 when $(P,A)
= \dagSix{•}{•}{•}{•}{•}{•}$ is essentially the same as the order
topology $\Opens_A(P)$ that we saw in sec.12 of \cite{OchsPH1}! This
holds for all graphs, and when $(P,A)$ is a 2CG --- for example, when
$$ (P,A) \;\; = \;\; \tcg{(P,A) 33} $$
we can abbreviate the result further using the ideas is sec.15
of \cite{OchsPH1}:
%
% (ph1p 28 "topologies-on-2CGs")
% (ph1     "topologies-on-2CGs")
%
$$ \Sub(1)
   \;\; = \;\;
   \Opens \tcg{(P,A) 33}
   \;\; = \;\;
   \zha{H 33}
$$

So: the ``logic'' of a topos of the form $\Set^{(P,A)}$ --- i.e., its
Heyting Algebra of subobjects of the terminal --- is exactly the
topology $\Opens_A(P)$.

\newpage

%  ____       _     /\    ______     _  __                         _         
% / ___|  ___| |_  |/\|  / /  _ \   / \ \ \   _ __ ___  _ __ _ __ | |__  ___ 
% \___ \ / _ \ __|      | || |_) | / _ \ | | | '_ ` _ \| '__| '_ \| '_ \/ __|
%  ___) |  __/ |_       | ||  __/ / ___ \| | | | | | | | |  | |_) | | | \__ \
% |____/ \___|\__|      | ||_| ( )_/   \_\ | |_| |_| |_|_|  | .__/|_| |_|___/
%                        \_\   |/       /_/                 |_|              
%
% «Set-PA-morphisms»  (to ".Set-PA-morphisms")
% \subsection{Morphisms in a topos $\Set^{(P,A)}$}
\subsection{Morphisms as natural transformations}

If $F$ and $G$ are objects of a category $\Set^\catA$ and $T:F→G$ is a
morphism between them then $F$ and $G$ are functors and $T:F→G$ is a
natural transformation, and T has to obey a ``naturalness condition''
that says that for every morphism $v:B→C$ in $\catA$ a certain
``obvious'' square must commute. We can draw that condition as the
commutativity of the middle square below,
%
%D diagram sqcond-1
%D 2Dx     100 +30   +30  +30    +45
%D 2D  100 B   FB -> GB   x |--> rx
%D 2D      |   |     |    -      -
%D 2D      |   |     |    |      v
%D 2D  +22 v   v     v    v      drx
%D 2D  +8  C   FC -> GC   dx |-> rdx
%D 2D
%D 2D  +20     F --> G
%D 2D
%D ren rx drx ==> (TB)(x) (Gv∘TB)(x)
%D ren dx rdx ==> (Fv)(x) (TC∘Fv)(x)
%D
%D (( B C -> .plabel= l v
%D    F G -> .plabel= a T
%D
%D    FB GB -> .plabel= a TB
%D    FB FC -> .plabel= l Fv
%D    GB GC -> .plabel= r Gv
%D    FC GC -> .plabel= a TC
%D
%D    x rx |-> rx drx |->
%D    x dx |-> dx rdx |->
%D ))
%D enddiagram
%D
$$\pu
  \diag{sqcond-1}
$$
and as the domain of $F$ and $G$ is $\Set$ we can express that
naturality as
$$∀(v:B→C). \, ∀x∈FB. \, (Gv∘TB)(x) = (TC∘Fv)(x)$$
and represent that as the square at the right above.

\msk

We will often draw these morphisms/natural transformations like this,
$$\pBargs{F_1}{F_2}{F_3}{F_4}{F_5}{F_6}
  \diagxyto/->/<250>^{T}
  \pBargs{F_1}{F_2}{F_3}{F_4}{F_5}{F_6}
$$
leaving the category $\catA$ implicit. The `$\diagxyto/->/^{T}$' is a
pack of six functions between sets, $T_1:F_1→G_1$, $\ldots$,
$T_6:F_6→G_6$ --- compare with the meaning of the `$\diagxyto/|->/$'
in sec.\ref{Set-PA}.

\msk

The definition of the local operator $j:Ω→Ω$ in
Figure \ref{fig:classifier-big} is a natural transformation written in
a very compact form. In that example $j_{3▁}(21) = 32$.

\directlua{tf_pop()}

% Local Variables:
% coding: utf-8-unix
% ee-tla: "joe"
% End:
   % (find-LATEX "2019J-ops-categories.tex")
\newpage
% (find-LATEX "2019J-ops-classifier.tex")
% (defun c () (interactive) (find-LATEXsh "lualatex -record 2019J-ops-classifier.tex" :end))
% (defun d () (interactive) (find-pdf-page "~/LATEX/2019J-ops-classifier.pdf"))
% (defun e () (interactive) (find-LATEX "2019J-ops-classifier.tex"))
% (defun u () (interactive) (find-latex-upload-links "2019J-ops-classifier"))
% (find-pdf-page   "~/LATEX/2019J-ops-classifier.pdf")
% (find-sh0 "cp -v  ~/LATEX/2019J-ops-classifier.pdf /tmp/")
% (find-sh0 "cp -v  ~/LATEX/2019J-ops-classifier.pdf /tmp/pen/")
%   file:///home/edrx/LATEX/2019J-ops-classifier.pdf
%               file:///tmp/2019J-ops-classifier.pdf
%           file:///tmp/pen/2019J-ops-classifier.pdf
% http://angg.twu.net/LATEX/2019J-ops-classifier.pdf
% (find-LATEX "2019.mk")

\directlua{tf_push("2019J-ops-classifier.tex")}

% «.Omega-and-j»		(to "Omega-and-j")
% «.pullbacks-formally»		(to "pullbacks-formally")
%   «.NTs-B-1»			(to "NTs-B-1")
%   «.NTs-B-C»			(to "NTs-B-C")
%   «.NTs-1-Om»			(to "NTs-1-Om")
%   «.NTs-C-Om»			(to "NTs-C-Om")
%   «.NTs-Om-Om»		(to "NTs-Om-Om")
%   «.fig:five-sqconds»		(to "fig:five-sqconds")
% «.pullbacks-visually»		(to "pullbacks-visually")

%   ___                              _     _ 
%  / _ \ _ __ ___     __ _ _ __   __| |   (_)
% | | | | '_ ` _ \   / _` | '_ \ / _` |   | |
% | |_| | | | | | | | (_| | | | | (_| |   | |
%  \___/|_| |_| |_|  \__,_|_| |_|\__,_|  _/ |
%                                       |__/ 
%
% «Omega-and-j»  (to ".Omega-and-j")
% (jopp 25 "Omega-and-j")
% (joe     "Omega-and-j")

\subsection{The classifier and the local operator}

We know that every category $\Set^{(P,A)}$ is a topos, but how do we
calculate and visualize its classifier object $Ω$ and the map $⊤:1→Ω$?
And what is the local operator $j:Ω→Ω$ ``associated to'' our
J-operator $J:\Sub(1)→\Sub(1)$?

\msk

%D diagram Omega-and-j
%D 2Dx     100   +30     +30
%D 2D  100 B --> 1
%D 2D      |     |
%D 2D      v     v
%D 2D  +30 C --> Om1 --> Om2
%D 2D
%D ren Om1 Om2 ==> Ω Ω
%D
%D (( B 1  -> .plabel= a !
%D    B C >-> .plabel= l i
%D    1 Om1 >-> .plabel= r ⊤
%D    C Om1 -> .plabel= a χ_B
%D    Om1 Om2 -> .plabel= a j
%D    B relplace 7 7 \pbsymbol{7}
%D ))
%D enddiagram
%D
%D diagram Omega-and-j-2
%D 2Dx     100   +30     +30
%D 2D  100 B ----------> 1
%D 2D      |             |
%D 2D      v             v
%D 2D  +30 C --> Om1 --> Om2
%D 2D
%D ren Om1 Om2 B ==> Ω Ω \ovl{B}
%D
%D (( B 1  -> .plabel= a !
%D    B C >-> .plabel= l \ovl{i}
%D    1 Om2 >-> .plabel= r ⊤
%D    C Om1 -> .plabel= a χ_B
%D    Om1 Om2 -> .plabel= a j
%D    B relplace 7 7 \pbsymbol{7}
%D ))
%D enddiagram
%D
\pu

% TODO: Explain the prequisites for this section. Explain that I
% learned this from Bell but McLarty is more readable.
%
% (find-books "__cats/__cats.el" "mclarty")
% (find-books "__cats/__cats.el" "bell")

We need to start by understanding two pullbacks. Remember that:

\begin{itemize}

\item $⊤:1→Ω$ has a property can be expressed in two equivalent ways:
  1) for each object $C$ we have $\Sub(C) ≅ \Hom(C,Ω)$, and 2) for
  every monic $B \monicto C$ there is exactly one map $χ_B:C→Ω$ making
  the square below --- ``the Q-shaped diagram'' --- a pullback:
  $$\diag{Omega-and-j}
  $$

\item a local operator (also called a
  ``modality'', a ``Lawvere-Tierney topology'', or a ``topology'') is
  a map $j:Ω→Ω$ obeying $j∘⊤=⊤$, $j∘j=j$ and $j∘∧=∧∘(j×j)$,

  % (find-books "__cats/__cats.el" "mclarty")
  % (find-mclartypage (+ 4 196) "21. Topologies")

\item a local operator $j$ induces a $j$-closure operator --- see
  chapter 21 of \cite{McLarty} or chapter 5 of \cite{BellLST} ---, and
  this $j$-closure operator can be seen as a map from each $\Sub(C)$
  to itself. The closure of a subobject $i: B \monicto C$ is the
  subobject $\ovl 1 : \ovl B \monicto C$ obtained by pullback in the
  diagram below (``the rectangle''):
  $$\diag{Omega-and-j-2}
  $$

\end{itemize}

We will write the restriction of a local operator $j$ to $\Sub(1)$ as
$\sfJ(j)$ and we will say that a $j$ is ``associated to'' a $J$ when
$\sfJ(j) = J$.

\msk

There are two ways to ``understand'' the pullbacks above: the first
one is by doing the calculations formally and checking that everything
works, the second one is by checking some particular cases and
developing visual intuition from that.

% In the next sections I will refer to the two diagrams above as ``the
% Q-shaped diagram'' and ``the rectangle''.

\newpage

%  ____  ____         __                            _ _       
% |  _ \| __ ) ___   / _| ___  _ __ _ __ ___   __ _| | |_   _ 
% | |_) |  _ \/ __| | |_ / _ \| '__| '_ ` _ \ / _` | | | | | |
% |  __/| |_) \__ \ |  _| (_) | |  | | | | | | (_| | | | |_| |
% |_|   |____/|___/ |_|  \___/|_|  |_| |_| |_|\__,_|_|_|\__, |
%                                                       |___/ 
%
% «pullbacks-formally»  (to ".pullbacks-formally")
\subsection{Understanding the pullbacks formally}
\label{pullbacks-formally}

The calculations are routine if we know the right language, and if we
suppose --- without loss of generality --- that the monix $i:B\monicto
C$ is a ``canonical subobject'' in the sense that each $B(p)⊆C(p)$ and
each function $B(p\ton!q):B(p)→B(q)$ is a restriction of the
corresponding function $C(p\ton!q):C(p)→C(q)$.

We need some definitions:
$$\begin{array}{rcl}
  1(p)          &=& \{*\} \\
  1(p\ton!q)(*) &=& * \\
  Ω(p)          &=& \Sub(↓p) \\
  Ω(p\ton!q)(R) &=& R∧↓q   \\%
  [5pt]
  ⊤(p)(*)       &=& ↓p     \\
  j(p)(R)       &=& R^*∧↓p \\
  χ_B(p)(R)     &=& \setofst{r∈↓p}{C(p\ton!r)(c)∈B(r)} \\
  \end{array}
$$

The first step is to check the five naturality conditions in the next
page --- we leave the rest to the reader. The main exercise is to
check that if the monic $i:B\monicto C$ is $i:P\monicto 1$ for a
truth-value $P$ then its closure is $i:\ovl P\monicto 1$ with $\ovl P$
being exactly $J(P)$, i.e., $P^*$.

\newpage

\pu
\begin{figure}[h!]
  \centering
  $\pu
   \scalebox{0.9}{$
   \begin{array}{l}
   \diag{B->1}   \\ \\
   \diag{B->C}   \\ \\
   \diag{1->Om}  \\ \\
   \diag{C->Om}  \\ \\
   \diag{Om->Om} \\
   \end{array}
   $}
  $
  \caption{The five square conditions in the Q-shaped diagram}
  \label{fig:five-sqconds}
\end{figure}

\newpage

%  ____  ____              _                 _ _       
% |  _ \| __ ) ___  __   _(_)___ _   _  __ _| | |_   _ 
% | |_) |  _ \/ __| \ \ / / / __| | | |/ _` | | | | | |
% |  __/| |_) \__ \  \ V /| \__ \ |_| | (_| | | | |_| |
% |_|   |____/|___/   \_/ |_|___/\__,_|\__,_|_|_|\__, |
%                                                |___/ 
%
% «pullbacks-visually»  (to ".pullbacks-visually")
\subsection{Understanding the pullbacks visually}
\label{pullbacks-visually}

The best way to develop visual intuition about the $Ω$ and the $j$
associated to a $((P,A),Q)$ is to try to work out the details in some
particular cases --- I've chosen two, presented as execises below.
They both use the $((P,A),Q)$, the $Ω$ and the $j$ from
Figure \ref{fig:classifier-big}.

\msk

{\bf Exercise 1.} In the case
%
%D diagram Omega-and-j-exercise-1-Q
%D 2Dx     100   +30     +30
%D 2D  100 B --> 1
%D 2D      |     |
%D 2D      v     v
%D 2D  +30 C --> Om1 --> Om2
%D 2D
%D ren B C ==> 11 33
%D ren Om1 Om2 ==> Ω Ω
%D
%D (( B 1  -> .plabel= a !
%D    B C >-> .plabel= l i
%D    1 Om1 >-> .plabel= r ⊤
%D    C Om1 -> .plabel= a χ_B
%D    Om1 Om2 -> .plabel= a j
%D    B relplace 7 7 \pbsymbol{7}
%D ))
%D enddiagram
%D
%D diagram Omega-and-j-exercise-1-rect
%D 2Dx     100   +30     +30
%D 2D  100 B ----------> 1
%D 2D      |             |
%D 2D      v             v
%D 2D  +30 C --> Om1 --> Om2
%D 2D
%D ren B C ==> \ovl{11} 33
%D ren Om1 Om2 ==> Ω Ω
%D
%D (( B 1  -> .plabel= a !
%D    B C >-> .plabel= l \ovl{i}
%D    1 Om2 >-> .plabel= r ⊤
%D    C Om1 -> .plabel= a χ_B
%D    Om1 Om2 -> .plabel= a j
%D    B relplace 7 7 \pbsymbol{7}
%D ))
%D enddiagram
%D
\pu
$$
\diag{Omega-and-j-exercise-1-Q}
\qquad
\diag{Omega-and-j-exercise-1-rect}
$$
what is $χ_B$? And what is $\ovl{11}$?

\msk

{\bf Exercise 2.} In the case
%
%D diagram Omega-and-j-exercise-2-Q
%D 2Dx     100   +30     +30
%D 2D  100 B --> 1
%D 2D      |     |
%D 2D      v     v
%D 2D  +30 C --> Om1 --> Om2
%D 2D
%D ren B C ==> 11 23
%D ren Om1 Om2 ==> Ω Ω
%D
%D (( B 1  -> .plabel= a !
%D    B C >-> .plabel= l i
%D    1 Om1 >-> .plabel= r ⊤
%D    C Om1 -> .plabel= a χ_B
%D    Om1 Om2 -> .plabel= a j
%D    B relplace 7 7 \pbsymbol{7}
%D ))
%D enddiagram
%D
%D diagram Omega-and-j-exercise-2-rect
%D 2Dx     100   +30     +30
%D 2D  100 B ----------> 1
%D 2D      |             |
%D 2D      v             v
%D 2D  +30 C --> Om1 --> Om2
%D 2D
%D ren B C ==> \ovl{11} 23
%D ren Om1 Om2 ==> Ω Ω
%D
%D (( B 1  -> .plabel= a !
%D    B C >-> .plabel= l \ovl{i}
%D    1 Om2 >-> .plabel= r ⊤
%D    C Om1 -> .plabel= a χ_B
%D    Om1 Om2 -> .plabel= a j
%D    B relplace 7 7 \pbsymbol{7}
%D ))
%D enddiagram
%D
\pu
$$
\diag{Omega-and-j-exercise-2-Q}
\qquad
\diag{Omega-and-j-exercise-2-rect}
$$

what is $χ_B$? And what is $\ovl{11}$?

% (elep 6 "elephant-A2.1.3")
% (ele    "elephant-A2.1.3")

% (elep 7 "elephant-A4.1.4")
% (ele    "elephant-A4.1.4")
% (elep 8 "elephant-A4.1.5")
% (ele    "elephant-A4.1.5")

% (ph1p 25 "topologies-as-partial-orders")
% (ph1     "topologies-as-partial-orders")

\directlua{tf_pop()}

% Local Variables:
% coding: utf-8-unix
% ee-tla: "joo"
% End:
   % (find-LATEX "2019J-ops-classifier.tex")
\newpage
% (find-LATEX "2019J-ops-kan.tex")
% (defun c () (interactive) (find-LATEXsh "lualatex -record 2019J-ops-kan.tex" :end))
% (defun d () (interactive) (find-pdf-page "~/LATEX/2019J-ops-kan.pdf"))
% (defun e () (interactive) (find-LATEX "2019J-ops-kan.tex"))
% (defun u () (interactive) (find-latex-upload-links "2019J-ops-kan"))
% (find-pdf-page   "~/LATEX/2019J-ops-kan.pdf")
% (find-sh0 "cp -v  ~/LATEX/2019J-ops-kan.pdf /tmp/")
% (find-sh0 "cp -v  ~/LATEX/2019J-ops-kan.pdf /tmp/pen/")
%   file:///home/edrx/LATEX/2019J-ops-kan.pdf
%               file:///tmp/2019J-ops-kan.pdf
%           file:///tmp/pen/2019J-ops-kan.pdf
% http://angg.twu.net/LATEX/2019J-ops-kan.pdf
% (find-LATEX "2019.mk")

% «.kan-extensions»	(to "kan-extensions")

\directlua{tf_push("2019J-ops-kan.tex")}

%  _  __                       _       
% | |/ /__ _ _ __     _____  _| |_ ___ 
% | ' // _` | '_ \   / _ \ \/ / __/ __|
% | . \ (_| | | | | |  __/>  <| |_\__ \
% |_|\_\__,_|_| |_|  \___/_/\_\\__|___/
%                                      
% «kan-extensions»  (to ".kan-extensions")
\subsection{Kan extensions}
\label     {kan-extensions}

\def\Lan{\text{Lan}}
\def\Ran{\text{Ran}}
\def\sfC{\mathsf{C}}
\def\sfD{\mathsf{D}}
\def\sfE{\mathsf{E}}

% (find-books "__cats/__cats.el" "riehl")
% (find-riehlccpage (+ 18  44) "1.7. The 2-category of categories")
% (find-riehlcctext (+ 18  44) "1.7. The 2-category of categories")
% (find-riehlccpage (+ 18  45) "Lemma 1.7.4 (horizontal composition)")
% (find-riehlcctext (+ 18  45) "Lemma 1.7.4 (horizontal composition)")
% (find-riehlccpage (+ 18  46) "whiskering")
% (find-riehlcctext (+ 18  46) "whiskering")
% (find-riehlccpage (+ 18 189) "6. All Concepts are Kan Extensions")
% (find-riehlccpage (+ 18 190) "6.1. Kan extensions")
% (find-riehlccpage (+ 18 190) "Dually, a right Kan")
% (find-riehlcctext (+ 18 190) "Dually, a right Kan")

In \cite{Riehl}, sec.6.1, right Kan extensions are explained using the
two diagrams below. The notation of cells is explained in sec.1.7 of
the book, and modulo the types --- that can be inferred from the
diagrams --- a right Kan extension of $K$ along $K$ is a pair $(\Ran_K
F,ε)$ such that for all $(G,α)$ there is a unique $β$ making
everything commute.
%
%D diagram riehl-ran-1
%D 2Dx     100 +40 +40
%D 2D  100 A0 ---> A2
%D 2D        ->  ->
%D 2D  +40     A1
%D 2D
%D ren A0 A1 A2 ==> \mathsf{C} \mathsf{D} \mathsf{E}
%D
%D (( A0 A2 -> .plabel= a F
%D    A0 A1 -> .plabel= l K
%D    A1 A2 -> .plabel= r G    .curve= _25pt
%D    A1 A2 varrownodes nil 17 nil <= .slide= -5pt .plabel= r δ
%D ))
%D enddiagram
%D
%D diagram riehl-ran-factored
%D 2Dx     100 +40 +40
%D 2D  100 A0 ---> A2
%D 2D        ->  -> ^
%D 2D  +40     A1 -/
%D 2D
%D ren A0 A1 A2 ==> \mathsf{C} \mathsf{D} \mathsf{E}
%D
%D (( A0 A2 -> .plabel= a F
%D    A0 A1 -> .plabel= l K
%D    A1 A2 -> .plabel= m \Ran_KF
%D    A1 A2 -> .plabel= r G    .curve= _25pt
%D    A0 A2 varrownodes 35 17 nil <=              .plabel= l ε
%D    A1 A2 varrownodes 20 17 nil <= .slide=  5pt .plabel= r β
%D ))
%D enddiagram
%D
$$\pu
  \diag{riehl-ran-1}
  \quad
  \diag{riehl-ran-factored}
$$

If we specialize $\sfE$ to $\Set$ and do some renamings, the diagram
becomes:
%
%D diagram my-ran-1
%D 2Dx     100 +40 +40
%D 2D  100 A0      A2
%D 2D
%D 2D  +40     A1
%D 2D
%D ren A0 A1 A2 ==> \catA \catB \Set
%D
%D (( A0 A2 -> .plabel= a D
%D    A0 A1 -> .plabel= l f
%D    A1 A2 -> .plabel= r C    .curve= _25pt
%D    A1 A2 varrownodes nil 17 nil <= .slide= -5pt .plabel= r α
%D ))
%D enddiagram
%D
%D diagram my-ran-2
%D 2Dx     100 +40 +40
%D 2D  100 A0      A2
%D 2D
%D 2D  +40     A1
%D 2D
%D ren A0 A1 A2 ==> \catA \catB \Set
%D
%D (( A0 A2 -> .plabel= a D
%D    A0 A1 -> .plabel= l f
%D    A1 A2 -> .plabel= m \Ran_fD
%D    A1 A2 -> .plabel= r C    .curve= _25pt
%D    A0 A2 varrownodes 35 17 nil <=              .plabel= l ε
%D    A1 A2 varrownodes 20 17 nil <= .slide=  5pt .plabel= r β
%D ))
%D enddiagram
%D
$$\pu
  \diag{my-ran-1}
  \quad
  \diag{my-ran-2}
$$
and if we change its {\sl shape} to stress that $ε$ ``looks like'' a
counit map and $\Ran_f$ ``looks like'' the right adjoint to the
functor $f^*$, we get this:
%
%D diagram geo-morph
%D 2Dx     100     +30 +35   +30     
%D 2D  100 L0      C0  C1    R1
%D 2D                          
%D 2D  +35 L2      C2  C3    R3
%D 2D                          
%D 2D  +20         C4  C5      
%D 2D
%D 2D  +20         C6  C7
%D 2D
%D ren C0 C1 C2 C3 C4 C5 ==> f^*C C D \Ran_fD \Set^\catA  \Set^\catB
%D ren C6 C7 ==> \catA \catB
%D ren L0 L2 ==> f^*\Ran_fD D
%D ren R1 R3 ==> C \Ran_ff^*C
%D
%D (( C0 C1 <-|
%D    C0 C2 -> .plabel= l \sm{β^\fl\\α}
%D    C1 C3 -> .plabel= r \sm{β\\α^♯}
%D    C2 C3 |->
%D    C0 C3 harrownodes nil 20 nil <-| sl^
%D    C0 C3 harrownodes nil 20 nil |-> sl_
%D
%D    C4 C5 <- sl^ .plabel= a f^*
%D    C4 C5 -> sl_ .plabel= b \Ran_f
%D
%D    C6 C7 -> .plabel= a f
%D    L0 L2 -> .plabel= l ε
%D    R1 R3 -> .plabel= r d
%D ))
%D enddiagram
%D
$$\pu
  \diag{geo-morph}
$$

When the categories $\catA$ and $\catB$ are finite posets we get:

\begin{itemize}

\item $\Set^\catA$ and $\Set^\catB$ are toposes (we saw this in
  sec.\ref{Set-PA}),
  %
  % (jopp 23 "Set-PA")
  % (joe     "Set-PA")

\item the functor $f^*$ is ``precomposition with $f$'', in this sense:
  if $C$ is an object of $\Set^B$ and $A∈\catA$ then $(f^*C)(A)$ is
  $C(f(A))$,

\item the left and right Kan extensions $\Lan_f$ and $\Ran_f$ and can be
  defined and calculated by the formulas in sec.6.2 of \cite{Riehl},
  %
  % (elep 7 "elephant-A4.1.4")
  % (ele    "elephant-A4.1.4")
  % (find-books "__cats/__cats.el" "riehl")
  % (find-riehlccpage (+ 18 193) "6.2. A formula for Kan extensions")

\item we have adjunctions $\Lan_f ⊣ f^* ⊣ \Ran_f$, and so the structure
  $(\Lan_f ⊣ f^* ⊣ \Ran_f)$ can be seen as an essential geometric
  morphism $f:\Set^\catA → \Set^\catB$ (\cite{EA}, A4.1.4); as $f^*$
  is a right adjoint it preserves limits (\cite{Riehl}, sec.4.5,
  and \cite{Awodey}, sec.9.6), and so $(f^* ⊣ \Ran_f)$ is a geometric
  morphism $f:\Set^\catA → \Set^\catB$. We usually rename $(\Lan_f ⊣
  f^* ⊣ \Ran_f)$ to $(f^! ⊣ f^* ⊣ f_*)$

\item when $f:\catA→\catB$ is something very simple we can find $\Ran_f D$
  ``by hand'' --- for example, in the example below, discussed
  in \cite{OchsACT2019}:
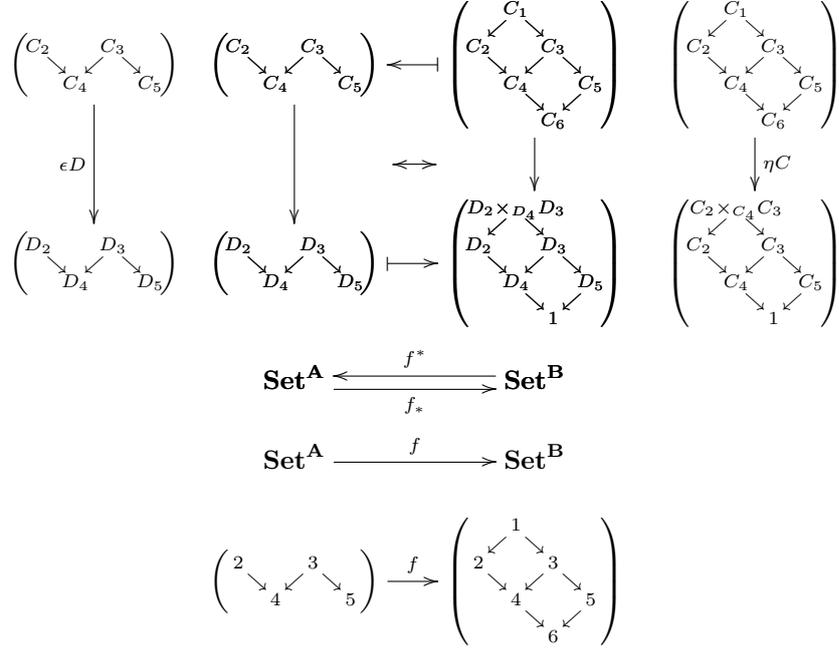
$$\pu
  \def\Ct{C_2 {×_{C_4}} C_3}
  \def\Dt{D_2 {×_{D_4}} D_3}
  \diag{internal-zgm-particular-case}
$$

\end{itemize}

\bsk

Every situation in which the category $\catB$ is a $(P,A)$ and the
category $\catA$ is the full subcategory of $(P,A)$ whose objects are
$P∖Q$ yields a situation like the one in the diagram above, in which
the maps $εD$ are isos, the geometric morphism $f$ is an ``inclusion''
and the functor that takes each $C$ to $f_*f^*C$ is a sheafification
functor. A diagram with an example fully worked out will be included
in the next version of this paper at the Arxiv.

% (find-books "__cats/__cats.el" "riehl")
% (find-riehlccpage (+ 18 136) "4.5. Adjunctions, limits, and colimits")
% (find-books "__cats/__cats.el" "awodey")
% (find-awodeyctpage (+ 10 197) "9.6 RAPL")

\newpage

\directlua{tf_pop()}

% Local Variables:
% coding: utf-8-unix
% ee-tla: "jok"
% End:
          % (find-LATEX "2019J-ops-kan.tex")

%L write_dnt_file()
\pu

% (find-LATEX "2019ilha-grande-poster-a4.tex" "references")

\newpage

\printbibliography

\end{document}

%  __  __       _        
% |  \/  | __ _| | _____ 
% | |\/| |/ _` | |/ / _ \
% | |  | | (_| |   <  __/
% |_|  |_|\__,_|_|\_\___|
%                        
% «make»  (to ".make")
% (find-LATEX "2019ilha-grande-poster-a4.tex")
% (find-LATEX "2019ilha-grande-poster-a4.tex" "write-poster-body")

 (eepitch-shell)
 (eepitch-kill)
 (eepitch-shell)

# Test the .zip.
# The e-script below downloads, unpacks and compiles the .zip in /tmp/edrx-latex/
#
rm -rfv /tmp/2019J-ops.zip
rm -rfv /tmp/edrx-latex/
cd /tmp/
# wget http://angg.twu.net/LATEX/2019J-ops.zip
cp -v ~/LATEX/2019J-ops.zip .
mkdir    /tmp/edrx-latex/
unzip -d /tmp/edrx-latex/ /tmp/2019J-ops.zip
cd       /tmp/edrx-latex/

cp -v ~/LATEX/2019J-ops-arxiv.tex 2019J-ops.tex

rm -v /tmp/2019J-ops-arxiv.zip
zip   /tmp/2019J-ops-arxiv.zip *

pdflatex 2019J-ops.tex
pdflatex 2019J-ops.tex

# (find-xpdfpage "/tmp/edrx-latex/2019J-ops.pdf")
# (find-fline    "/tmp/2019J-ops-arxiv.zip")

 (eepitch-shell)
 (eepitch-kill)
 (eepitch-shell)
# (find-LATEXfile "2017planar-has-1.mk")
make -f 2017planar-has-1.mk STEM=2019J-ops-arxiv veryclean
make -f 2017planar-has-1.mk STEM=2019J-ops-arxiv pdf
# (jopp)
# (jopp 25)

% Local Variables:
% coding: utf-8-unix
% ee-tla: "NONE"
% End: